\documentclass[10pt]{article}

\usepackage{geometry}
\usepackage{amssymb,amsmath} 
\usepackage{graphicx}        
\usepackage{float}  		
\usepackage{url}
\usepackage{float} 
\usepackage{psfrag} 
\usepackage{url}
\usepackage{subfig}

\usepackage{chemarr}
\usepackage{fancybox}
\usepackage{multicol}
\usepackage{wrapfig}
\usepackage{stfloats}
\usepackage{mathtools, cuted}
\usepackage{url}
\usepackage{cite}

\usepackage{epstopdf}
\usepackage{xcolor}
\usepackage{tabularx}
\usepackage{booktabs} 

\usepackage{soul}

\renewcommand\thesection{\Roman{section}}
\renewcommand\thesubsection{\thesection-\Alph{subsection}}

\begin{document}

\title{State observation and sensor selection for nonlinear networks}
\author{Aleksandar~Haber,
		Ferenc~Molnar,
        and~Adilson~E.~Motter
\footnote{This research was supported in part by MURI Grant No.~ARO-W911NF-14-1-0359, in part by  Simons Foundation Award No.~342906, and in part by NCI CR-PSOC Grant No.~1U54CA193419.
A. Haber was with the Department of Physics and Astronomy, Northwestern University, Evanston, IL 60208 USA, when this research was performed. He is now with the Department of Engineering Science and Physics, City University of New York, College of Staten Island, Staten Island, NY 10314 USA (e-mail: aleksandar.haber@csi.cuny.edu).
F. Molnar and A. E. Motter are with the Department of Physics and Astronomy, Northwestern University, Evanston, IL 60208 USA (e-mails: ferenc.molnar@northwestern.edu; motter@northwestern.edu).
}}

%


\maketitle

\begin{abstract}
A large variety of dynamical systems, such as chemical and biomolecular systems, can be seen as networks of nonlinear entities. Prediction, control, and identification of such nonlinear networks require knowledge of the state of the system. However, network states are usually unknown, and only a fraction of the state variables are directly measurable. The observability problem concerns reconstructing the network state from this limited information. Here, we propose a general optimization-based approach for observing the states of nonlinear networks and for optimally selecting the observed variables. Our results reveal several fundamental limitations in network observability, such as the trade-off between the fraction of observed variables and the observation length on one side, and the estimation error on the other side. We also show that, owing to the crucial role played by the dynamics, purely graph-theoretic observability approaches cannot provide conclusions about one's practical ability to estimate the states. We demonstrate the effectiveness of our methods by finding the key components in biological and combustion reaction networks from which we determine the full system state. Our results can lead to the design of novel sensing principles that can greatly advance prediction and control of the dynamics of such networks.
\end{abstract}


\clearpage

\section{Introduction}

Reaction systems, biophysical networks, and power grids are typical examples of systems with nonlinear network dynamics. 
The knowledge of the network state is important for the prediction \cite{camacho2013model},  
control \cite{khalil2002nonlinear,whalen2015observability,yan2015spectrum,lin2013design}, 
and identification \cite{timme2007revealing,verhaegen2007,haber2014subspace} of such systems.  
Determining the network state is challenging in practice because one is generally able to measure the time-series of only a fraction of all state variables; when the complete state can be determined from this limited information,  the network is said to be observable \cite{isidori2013nonlinear}. 
The problem of reconstructing the network state can be divided into two parts: ({\it i}\,) selection of state variables that need to be measured in order to guarantee the network observability; ({\it ii}\,) design of a state reconstructor (or observer) on the basis of the state variables selected in the first part. Despite the recent interest in the literature, problems ({\it i}) and ({\it ii}) remain open for nonlinear networks.

The classical approaches for the observability analysis of nonlinear systems rely on Lie-algebraic formulations \cite{isidori2013nonlinear}. However, these formulations cannot be used to optimally select the state variables (\textit{sensors}) guaranteeing network observability. On the other hand, the problem of selecting a minimal number of sensors that may guarantee structural observability of the network has been considered in \cite{liu2013observability}. Structural observability concerns the study of the connectivity between the state variables and outputs, without taking into account the precise values of the model parameters. In \cite{liu2013observability} the structural observability problem was considered by examining the observability inference diagram (OID),
which is a graph representing the dependences between the 
variables. The OID is constructed for network dynamics described by coupled first-order ordinary differential equations by choosing the state variables as nodes and adding a directed 
edge from node $i$ to node $j$ if variable $j$ appears on the r.h.s. of the equation for variable $i$. 
By analyzing the structure of this graph in terms of its strongly connected components (SCCs), it is possible to draw conclusions on the number and location of sensors to guarantee structural observability, namely that the minimal sets consist of one sensor in each root SCC of the OID (a root SCC is an SCC with no incoming edges). This approach offers an elegant 
graph-theoretic contribution to the structural observability problem. Purely structure-based approaches have also been proposed for the observation and reconstruction of attractor dynamics \cite{mochizuki2013,fielder2013}. These graphical approaches are successful in providing insights into the relation between network topology and observability. However, since these approaches do not explicitly take into account model parameters, they are not designed to guarantee near optimal performance of the state reconstruction.

The optimal selection of control nodes in networks with {\it linear} dynamics
has been extensively studied in the literature (see, e.g., \cite{summers2016submodularity}). Since that problem is dual to the problem of sensor selection (problem ({\it i}\,)), in linear networks, the methods from these previous studies can be used for sensor selection while accounting for the relevant model parameters. In \textit{nonlinear} networks, however,
optimal sensor selection still remains an open problem. For example, methods based on empirical Gramians in low-dimensional systems \cite{qi2015optimal,lall2002subspace,krener2009measures,powel2015empirical,singh2005determining} are not applicable to large-scale networks due to their high computational complexity and, as we show in this paper, low accuracy under realistic conditions. State estimation (problem  ({\it ii}\,)), on the other hand, has been studied extensively in nonlinear systems, and various approaches have been proposed, such as nonlinear extensions of the Kalman filter \cite{julier1997new,simon2006optimal}, particle filters \cite{stano2013saturated,stano2013parametric}, moving horizon estimation (MHE) techniques \cite{moraal1995observer,alessandri2008moving}, and others \cite{abarbanel2009}.
However, the applicability of such approaches to large-scale nonlinear networks  has not been investigated under the realistic conditions of a limited number of sensor nodes and a limited observation horizon.


Here, we propose a \textit{unified}, optimization-based framework for observing the states and optimally selecting the sensors in nonlinear networks, thereby offering a general solution to both problem ({\it i}\,) and problem ({\it ii}\,) \textit{under the same framework}.
We adopt the basic formulation of the open-loop MHE approach \cite{moraal1995observer,alessandri2008moving}, and formulate the state estimation problem as an optimization problem. Consequently, our approach can easily take into account various state constraints (e.g., min-max bounds and even nonlinear constraints). Moreover, the MHE approach enables us to study the influence of the observation horizon on the state estimation performance.

To the best of our knowledge, our approach is the most scalable procedure currently available
for sensor selection in nonlinear systems (not only in networks). We present extensive comparisons with existing approaches for validation. In addition, unlike other state-of-the-art methods for nonlinear state estimation \cite{simon2006optimal,rawlings2006,patwardhan2012}, our approach is capable of explicitly accounting for stiff nonlinear dynamics in a computationally efficient manner.

\begin{figure*}[h!]
\center
\includegraphics[width=\textwidth]{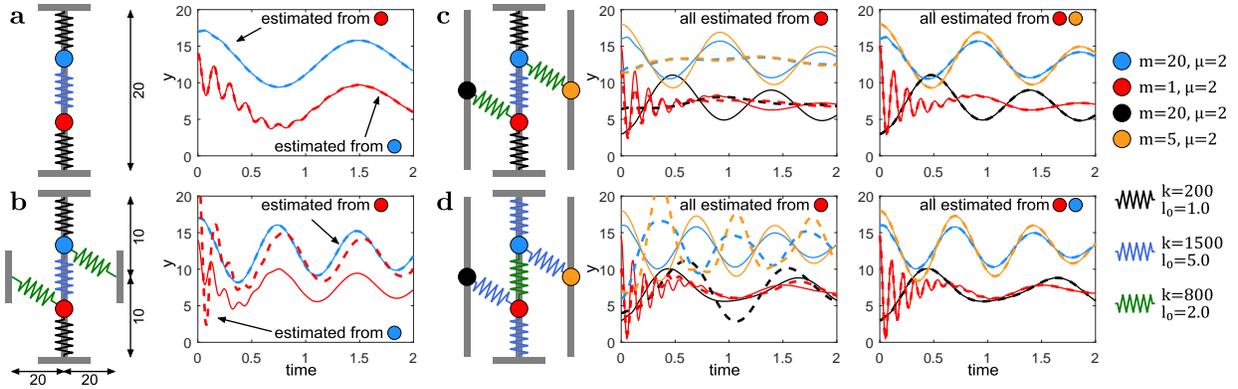}
\caption{\small State estimation in spring-mass networks.
Networks of (\textbf{a}) two masses subject to linear forces,  (\textbf{b}) two masses subject to nonlinear forces, and  (\textbf{c}, \textbf{d}) four masses subject to nonlinear forces. In each case, the point masses are restricted to move vertically and the forces (including the nonlinear ones) are emulated by linear springs in the plane. The  spring constants $k$, rest lengths $l_0$, friction coefficients $\mu$, and masses  $m$ are color-coded (legend on the right); the dimensions of the systems are marked on the figure (the dimensions in \textbf{c} and \textbf{d} are the same as in \textbf{b}). The initial states (position and velocity) of each mass are estimated only from the direct observation of the position of the subset of masses marked on the plots. The plots compare the true trajectories over time (solid lines) and those calculated from estimated initial states (dashed lines), color-coded as the masses. In \textbf{a} (linear case), estimation is successful from the observations of either mass, whereas in \textbf{b} (nonlinear case), estimation is  successful only if the smaller mass is observed. In \textbf{c} and \textbf{d} (larger networks), estimation is only successful if at least two masses are observed; comparison between \textbf{c} and \textbf{d} further shows that the optimal sensor placement (directly observed masses) depends not only on the OID but also on the dynamical parameters. Note that this is the case even though each network has as single (root) SCC.}
\label{fig1}
\end{figure*}

Our solution reveals the significant implication that, by virtue of realistic limitations in numerical and modeling precision, explicit state determination often requires a larger number of sensors than predicted from graph-theoretic approaches; moreover, the sensor nodes can depend strongly on the dynamical parameters even when the OID remains the same. This is illustrated in Fig.~\ref{fig1} for simple networks within the framework we present below. 

We validate our approach by performing extensive numerical experiments on biological
\cite{saadatpour2011dynamical,calzone2010mathematical, zhang2008network} 
and combustion reaction \cite{turns1996introduction,maas1992simplifying,smirnov2014modeling,perini2012analytical,o2004comprehensive,smith1999gri}
networks, which are examples par excellence of systems with nonlinear (and also stiff) dynamics. In particular, we specifically selected networks whose control is a subject of current research \cite{zanudo2015cellprog,cornelius2013realistic,wells2015control}. The 
numerical results enable us to detect which species concentrations and genes/gene products 
are the most important for the accurate state determination of these networks, thus demonstrating the efficacy of our method to reconstruct network states from limited measurement information.

The paper is organized as follows. In Section \ref{section2} we postulate the models we consider, and define the state estimation and sensor selection problems. Our approach to 
state estimation is presented in Section \ref{section3}, while the results on the optimal sensor selection are detailed in Section \ref{section4}. In Section \ref{section-results} we present and discuss our numerical experiments on the combustion and biological networks. Our
conclusions are presented in Section \ref{section-conclude}.

\section{Problem formulation}
\label{section2}

\begin{table}[b!]
\caption{Table of definitions}
\label{table1}
\centering
\begin{tabularx}{5in}{l X}
\toprule
\textbf{symbol} & \textbf{description}\\
\midrule
$\mathbf{x}(t)$ & vector of state variables at time $t$\\
$\mathbf{y}(t)$ & output vector at time $t$\\
$\mathbf{q}(\cdot)$ & vector function defining the dynamical equation of the system \\
$\text{col}\left(\mathbf{x}_{1},..., \mathbf{x}_{N}\right)$ & 
vector formed by collating the (column) vectors $\mathbf{x}_{1},..., \mathbf{x}_{N}$ \\
$n$ & number of nodes (state variables)
in the network \\
$r$ & number of directly observed nodes (state variables)\\
$N$ & observation length\\
$f$ & fraction of directly observed nodes (state variables) \\
$A=[a_{ij}]$ & matrix with entries $a_{ij}$ \\
$I_k$ & identity matrix of size $k \times k$\\
\bottomrule
\end{tabularx}
\end{table}

We focus on the general class of nonlinear networks described by a model of the form
\begin{align}
\dot{\mathbf{x}}(t)= \mathbf{q}\left(\mathbf{x}(t)\right),
\label{continiousTimeGeneral}
\end{align}
where $\mathbf{q}\left(\mathbf{x}\right):\mathbb{R}^{n} \rightarrow  \mathbb{R}^{n} $ is 
a nonlinear function at least twice continuously differentiable, and $\mathbf{x}\in \mathbb{R}^{n}$ represents the network state. 
Without loss of generality, we assume that a single state variable is associated to each node in the network.
With the model in equation~\eqref{continiousTimeGeneral} we associate a measurement equation:
\begin{align}
\mathbf{y}(t)=C\mathbf{x}(t),
\label{observationEquation}
\end{align}
where $\mathbf{y}(t)\in \mathbb{R}^{r}$, $r\le n$, is the output (measurement) vector and $C\in\mathbb{R}^{r\times n}$. The essential notation used in this paper is summarized in Table~\ref{table1}. The entries of $C$ are zero except for a single entry of $1$ in each row, corresponding to a sensor. 
For simplicity and brevity, in this paper we do not consider the effect of noise in equation~\eqref{observationEquation}; however this effect can be evaluated using existing methods.

We develop a state estimation method on the basis of models obtained by discretizing the continuous-time model in equation~\eqref{continiousTimeGeneral}. The discretization of the continuous dynamics at an early stage of the estimator design is a standard technique \cite{verhaegen2007,moraal1995observer,alessandri2008moving,simon2006optimal} that simplifies calculations, and is justified by the discreteness of the data usually obtained in real experiments \cite{verhaegen2007}. In particular, we postulate models based on discretization 
schemes that can perform well even if the continuous-time model in equation~\eqref{continiousTimeGeneral} is stiff, which is often the case, for example, for reaction networks. We consider models based on the backward Euler (BE), trapezoidal implicit (TI) and two-stage implicit Runge-Kutta (IRK) discretization techniques \cite{iserles2009first}. The BE technique leads to the model $\mathbf{x}_{k}=\mathbf{x}_{k-1}+h\mathbf{q}\left(\mathbf{x}_{k}\right)$, where $h>0$ is a discretization step, $\mathbf{x}_{k}=\mathbf{x}(kh)$, and $k=0,1,\ldots$, is a discrete-time step. The TI technique leads to the model $\mathbf{x}_{k}=\mathbf{x}_{k-1}+0.5h\left( \mathbf{q}\left(\mathbf{x}_{k}\right)+\mathbf{q}\left(\mathbf{x}_{k-1}\right) \right)$.
Finally, the model postulated on the basis of the IRK technique is
\begin{align}
\pmb{\zeta}_{1,k}&=\mathbf{x}_{k-1}+(h/4)\left(\mathbf{q}\left(\pmb{\zeta}_{1,k}\right)-\mathbf{q}\left(\pmb{\zeta}_{2,k}\right)\right) \label{rk1}, \\
\pmb{\zeta}_{2,k}&=\mathbf{x}_{k-1}+(h/12)\left(3\mathbf{q}\left(\pmb{\zeta}_{1,k}\right)+5\mathbf{q}\left(\pmb{\zeta}_{2,k}\right)\right) \label{rk2},\\ 
\mathbf{x}_{k}&=\mathbf{x}_{k-1}+(h/4)\left(\mathbf{q}\left(\pmb{\zeta}_{1,k}\right)+3\mathbf{q}\left(\pmb{\zeta}_{2,k}\right)\right) \label{rk3},
\end{align}
where $\pmb{\zeta}_{1,k},\pmb{\zeta}_{2,k}\in \mathbb{R}^{n}$ are vectors needed to compute $\mathbf{x}_{k}$ 
once vector $\mathbf{x}_{k-1}$ has been determined.
With the introduced models we associate an output equation,
\begin{align}
\mathbf{y}_{k}&=C\mathbf{x}_{k},\label{implicitCompactOutput}
\end{align} 
where $\mathbf{y}_{k}\in \mathbb{R}^{r}$ is defined analogously to $\mathbf{x}_{k}$. 

Models of real networks often involve modeling uncertainties. To emulate these uncertainties, we use the model in equation~\eqref{continiousTimeGeneral} as a data-generating model, acting as a ``real'' physical system whose state we want to estimate. The observation data 
are generated by numerically integrating equation~\eqref{continiousTimeGeneral} using a more accurate technique than the ones considered above. This way, we are able to validate the performance of our observation strategy against the model uncertainties originating from the difference between the discrete-time model formulations and the model in equation~\eqref{continiousTimeGeneral}.

Using the model formulations described above, our first goal is to estimate the initial state $\mathbf{x}_{0}$ from 
a set of observations $\mathcal{O}_{N}=\{\mathbf{y}_{0},\mathbf{y}_{1},\mathbf{y}_{2},\ldots, \mathbf{y}_{N-1} \}$, where $N$ is the observation length. The illustration in Fig.~\ref{fig1}, for example, was generated using the IRK model for $h=0.01$ and $N=200$.
Our second goal is to choose an optimal set of $r$ sensor nodes
that allows for the most accurate reconstruction of the initial state, where this set can be 
constrained not to include specific nodes in the network.

\section{Initial state estimation}
\label{section3}

To begin, from equation~\eqref{implicitCompactOutput} we define
\vspace{-1mm}
\begin{align}
\begin{array} {rl}
\mathbf{g}_{0} \; := & \mathbf{y}_{0}-C\mathbf{x}_{0}, \\
\mathbf{g}_{1} \; := & \mathbf{y}_{1}-C\mathbf{x}_{1}, \\
&  \hdots, \\
\mathbf{g}_{N-1} \; := & \mathbf{y}_{N-1}-C\mathbf{x}_{N-1}. 
\end{array}  
\label{eqN}
\end{align}
Using the column vector $\mathbf{g}=\text{col}\left(\mathbf{g}_{0},\mathbf{g}_{1},\ldots, \mathbf{g}_{N-1} \right)\in \mathbb{R}^{Nr}$, and  on the basis of equation~\eqref{eqN}, we can define the following equation:
\begin{align}
\mathbf{g}\left(\mathbf{x}_{0} \right)=0,
\label{nonlinearCompact} 
\end{align}
where $\mathbf{g}\left(\mathbf{x}_{0} \right):\mathbb{R}^{n}\rightarrow \mathbb{R}^{Nr}$ is a nonlinear vector function of the initial state 
(i.e., equation~\eqref{nonlinearCompact} represents a system of nonlinear equations). The vector $\mathbf{g}$ is only a function of $\mathbf{x}_{0}$ because the states in the sequence $\{\mathbf{x}_{1},\mathbf{x}_{2},\ldots,\mathbf{x}_{N-1} \}$ are coupled together through the postulated state-space models, and they depend only on $\mathbf{x}_{0}$. To proceed, it is beneficial to introduce the notation $\mathbf{y}=\text{col}\left(\mathbf{y}_{0},\mathbf{y}_{1},\ldots, \mathbf{y}_{N-1} \right)\in \mathbb{R}^{Nr}$ and  $\mathbf{w}=\text{col}\left(C\mathbf{x}_{0},C\mathbf{x}_{1},\ldots, C\mathbf{x}_{N-1} \right)\in \mathbb{R}^{Nr}$, where $\mathbf{w}=\mathbf{w}\left(\mathbf{x}_{0} \right)$. From equation~\eqref{nonlinearCompact}, it follows that
\begin{align}
\mathbf{y}&=\mathbf{w}\left(\mathbf{x}_{0}\right). \label{nonlinearCompact3}
\end{align}
The network is observable if the initial state $\mathbf{x}_{0}$ can be \textit{uniquely} determined from the set of observations $\mathcal{O}_{N}$. A formal definition of the observability of discrete-time systems is given in \cite{hanba2009uniform}: \textit{the system is said to be uniformly observable on a set $\Omega\subset \mathbb{R}^{n}$ if $\exists N>0$ such that the map $\mathbf{w}\left(\mathbf{x}_{0}\right)$ is injective as a function of $\mathbf{x}_{0}$.} 
A sufficient condition for observability in $\Omega$ for some $N$ is given by the rank condition $\mathsf{ rank}\left(J\left(\mathbf{x}_{0} \right)\right)=n\; \forall\, \mathbf{x}_{0} \in \Omega$, where $J\left(\mathbf{x}_{0} \right) \in \mathbb{R}^{Nr \times  n}$ is the Jacobian matrix of the map $\mathbf{w}\left(\mathbf{x}_{0} \right)$ \cite{hanba2009uniform}.
In order to satisfy the rank condition, the Jacobian matrix should be at least a square matrix, which tells us that the observation length $N$ should satisfy $Nr\ge n$.
Beyond justifying this constraint, however, the rank condition cannot be applied directly because $\mathbf{x}_{0}$ is unknown a priori.

It is immediate that the existence and uniqueness of the solution of equation~\eqref{nonlinearCompact} inside of a domain guarantee the observability in this domain. For square systems 
in equation~\eqref{nonlinearCompact} (i.e., $Nr=n$), the Kantorovich theorem \cite{dennis1996numerical} gives us a condition for the existence and the uniqueness of the solution. 
In particular, the Kantorovich theorem tells us that equation~\eqref{nonlinearCompact} has a unique solution in an Euclidean ball around $\mathbf{x}_{0}^{(0)}$, where this vector is an initial guess of the Newton method for solving the equation, if 
1) the Jacobian of $\mathbf{g}$ is nonsingular at $\mathbf{x}_{0}^{(0)}$, 
2) the Jacobian is Lipschitz continuous in a region containing the initial guess, and 
3) the first step of the Newton method taken from the initial point is relatively small. Because the Jacobians of $\mathbf{g}\left(\mathbf{x}_{0} \right)$ and $\mathbf{w}\left(\mathbf{x}_{0} \right)$ only differ by a sign, it follows
that for square systems the observability condition based on the rank of the Jacobian at $\mathbf{x}_{0}$ can be substituted by a rank condition at $\mathbf{x}_{0}^{(0)}$. 
Therefore, satisfying the rank condition at $\mathbf{x}_{0}^{(0)}$ (i.e., $\mathsf{ rank}(J(\mathbf{x}_{0}^{(0)} ))=n$) implies that equation~\eqref{nonlinearCompact} has a unique solution around $\mathbf{x}_{0}^{(0)}$. 

We determine the solution of equation~\eqref{nonlinearCompact} by numerically solving the following optimization problem:
\begin{align}
& \min_{\mathbf{x}_{0}} \left\|\mathbf{g}(\mathbf{x}_{0})  \right\|_{2}^{2}, \;\; \text{subject to}\;\; \underline{\mathbf{x}}_{0} \le \mathbf{x}_{0} \le \overline{\mathbf{x}}_{0},
\label{nonlinearOptimization} 
\end{align}
where $\le$ is applied element-wise and $\underline{\mathbf{x}}_{0},\overline{\mathbf{x}}_{0} \in \mathbb{R}^{n}$ are the bounds on the optimization variables taking into account the physical constraints of the network (these constraints can be modified to also include nonlinear functions of $\mathbf{x}_{0}$).
The problem in equation~\eqref{nonlinearOptimization} is a constrained non-linear least-squares problem \cite{dennis1996numerical}, which we solve using the trust region reflective (TRR) algorithm \cite{byrd1988approximate,more1983computing,coleman2001preconditioned,sorensen1997minimization}.
It is important to note that while we do not consider measurement noise here, the effect of the noise on least-squares problems is well-studied in the literature \cite{ljung1999system}.
To quantify the estimation accuracy, we introduce the estimation error $\eta=\left\|\hat{\mathbf{x}}_{0}-\mathbf{x}_{0} \right\|_{2}/\left\|\mathbf{x}_{0}   \right\|_{2}$, where $\hat{\mathbf{x}}_{0}$ is the solution of the optimization problem~\eqref{nonlinearOptimization}.

\subsection{Calculating the Jacobian}

In order to significantly speed up the computations of the state estimate, we derive the Jacobians of the function $\mathbf{g}(\mathbf{x}_{0})$
in equation~\eqref{nonlinearOptimization}  for the three postulated models. 
Without these analytical formulas, each element of the Jacobian would have to be computed using finite differences, resulting in unwanted computational burden.
The Jacobian matrix $J(\mathbf{x}_{0})\in \mathbb{R}^{Nr \times n}$ of the function $\mathbf{g}\left(\mathbf{x}_{0}\right)\in \mathbb{R}^{Nr}$ 
needs to be calculated at the point $\mathbf{x}_{0}^{(i)}$, where $i$ indicates the $i$th iteration of the TRR method for solving the optimization problem in equation~\eqref{nonlinearOptimization}. 
For the computation of the Jacobians we adopt the numerator layout notation.
Using the chain rule, it is possible to express the Jacobian matrix $J(\mathbf{x}_{0}^{(i)})$ as follows:
\begin{align}
J(\mathbf{x}_{0}^{(i)}) = \begin{bmatrix} -C \\ -C\frac{\partial \mathbf{x}_{1}}{\partial\mathbf{x}_{0}} \\ -C\frac{\partial \mathbf{x}_{2}}{\partial \mathbf{x}_{1}}\frac{\partial \mathbf{x}_{1}}{\partial\mathbf{x}_{0}} \\ \vdots \\
-C\frac{\partial \mathbf{x}_{N-1}}{\partial \mathbf{x}_{N-2}}\frac{\partial \mathbf{x}_{N-2}}{\partial \mathbf{x}_{N-3}}\ldots \frac{\partial \mathbf{x}_{1}}{\partial \mathbf{x}_{0}}  \end{bmatrix},
\label{jacobianOriginal}
\end{align}
 where
\begin{align}
\frac{\partial \mathbf{x}_{j+1}}{\partial\mathbf{x}_{j}}=\frac{\partial \mathbf{x}_{j+1}}{\partial\mathbf{x}_{j}}\Bigr|_{ \mathbf{x}_{j}^{(i)}}.
\label{derivativesForm2}
\end{align}
That is, the derivatives of the form of $\partial \mathbf{x}_{j+1} /  \partial  \mathbf{x}_{j}$, $j=0,1,\ldots, N-1$, in equation \eqref{jacobianOriginal} are evaluated at $\mathbf{x}_{j}^{(i)}$. Note that in equation \eqref{jacobianOriginal} we assume that $\mathbf{x}_{j+1}$ is formally expressed as a function of only $\mathbf{x}_{j}$ (without any implicit dependence on $\mathbf{x}_{j+1}$). Under mild conditions, the existence of such an expression follows from the implicit function theorem.

\subsubsection*{Jacobians for the BE and TI models}
The first challenge in computing the Jacobians originates from the implicit nature of the state equations of the BE and TI models defined in Section~\ref{section2}. Namely, the term $\mathbf{x}_{k}$ appears on both sides of the corresponding equations, and consequently, the corresponding partial derivatives will appear on both sides of equations. To see this in the case of the BE model, 
for the time step $j$ we can write
\begin{align}
\frac{\partial \mathbf{x}_{j}}{\partial\mathbf{x}_{j-1}}=I_n+h\frac{\partial \mathbf{q}(\mathbf{x}_{j})}{\partial \mathbf{x}_{j}}\Bigr|_{ \mathbf{x}_{j}^{(i)}}\frac{\partial \mathbf{x}_{j}}{\partial \mathbf{x}_{j-1}},
\label{JacobianEvaluation1}
\end{align}
where $j=1,2,\ldots,N-2$ 
and $I_n$ is the $n\times n$ identity matrix.
From this equation it follows that 
\begin{align}
\underbrace{\left(I_n-h\frac{\partial \mathbf{q}(\mathbf{x}_{j})}{\partial \mathbf{x}_{j}}\Bigr|_{ \mathbf{x}_{j}^{(i)}} \right)}_{A_{1}}\frac{\partial \mathbf{x}_{j}}{\partial\mathbf{x}_{j-1}}=I_n,
\label{JacobianEvaluation2}
\end{align}
where $A_{1}\in \mathbb{R}^{n\times n}$. Assuming that the matrix $A_{1}$ is invertible (which is a sufficient condition for the implicit function theorem to guarantee the existence of $\mathbf{x}_{j}$ as a function of only $\mathbf{x}_{j-1}$), we obtain
\begin{align}
\frac{\partial \mathbf{x}_{j}}{\partial\mathbf{x}_{j-1}}=A_{1}^{-1}.
\label{JacobianEvaluation3}
\end{align}
Using a similar argument, in the case of the TI model, we obtain
\begin{equation}
\frac{\partial \mathbf{x}_{j}}{\partial\mathbf{x}_{j-1}}=\Big(I_n-0.5h\frac{\partial \mathbf{q}(\mathbf{x}_{j})}{\partial \mathbf{x}_{j}}\Bigr|_{\mathbf{x}_{j}^{(i)}} \Big)^{-1}\Big(I_n+0.5h\frac{\partial \mathbf{q}(\mathbf{x}_{j-1})}{\partial \mathbf{x}_{j-1}}\Bigr|_{\mathbf{x}_{j-1}^{(i)}} \Big).
\label{JacobianEvaluation4}
\end{equation}
The expressions in equations~\eqref{JacobianEvaluation3} and \eqref{JacobianEvaluation4} show that to evaluate the Jacobian matrices of the corresponding models at the point $\mathbf{x}_{0}^{(i)}$, one actually needs to know the values of $\mathbf{x}_{1}^{(i)},\ldots,\mathbf{x}_{N-1}^{(i)}$. These values can be obtained by simulating the BE model or the TI model, starting from the initial point $\mathbf{x}_{0}^{(i)}$. This procedure needs to be repeated for every iteration $i$ of the TRR method to solve the optimization problem in equation~\eqref{nonlinearOptimization}. 


\subsubsection*{Jacobian for the IRK model}
In the case of the IRK model, defined in equations~\eqref{rk1}, \eqref{rk2}, and \eqref{rk3}, 
the evaluation of the Jacobian matrix becomes even more involved numerically. By setting the time step $k$ in equation~\eqref{rk3} 
to $j$, and differentiating such an expression, we obtain
\begin{equation}
\frac{\partial \mathbf{x}_{j}}{\partial\mathbf{x}_{j-1}} = I_n+\frac{h}{4}
 \frac{\partial\mathbf{q}( \pmb{\zeta}_{1,j})}{\partial \pmb{\zeta}_{1,j}}\Bigr|_{\pmb{\zeta}_{1,j}^{(i)}}\frac{\partial \pmb{\zeta}_{1,j}}{\partial \mathbf{x}_{j-1}}+\frac{3h}{4}\frac{\partial \mathbf{q}( \pmb{\zeta}_{2,j})}{\partial \pmb{\zeta}_{2,j}}\Bigr|_{\pmb{\zeta}_{2,j}^{(i)}}\frac{\partial \pmb{\zeta}_{2,j}}{\partial \mathbf{x}_{j-1}}
\label{trapezoidalPD1}
\end{equation}
To evaluate \eqref{trapezoidalPD1} we need to determine the partial derivatives $\partial \pmb{\zeta}_{1,j}/\partial \mathbf{x}_{j-1}$ and $\partial \pmb{\zeta}_{2,j}/\partial \mathbf{x}_{j-1}$. By differentiating equations~(3) and (4), 
we obtain
\begin{align}
\underbrace{\begin{bmatrix}\frac{\partial \pmb{\zeta}_{1,j}}{\partial \mathbf{x}_{j-1}}\\ \frac{\partial \pmb{\zeta}_{2,j}}{\partial \mathbf{x}_{j-1}} \end{bmatrix}}_{S} &=\underbrace{\begin{bmatrix}I_n  \\  I_n  \end{bmatrix}}_{\widetilde{I}_n}+\underbrace{\begin{bmatrix} \frac{h}{4} \frac{\partial\mathbf{q}( \pmb{\zeta}_{1,j})}{\partial \pmb{\zeta}_{1,j}} &  -\frac{h}{4}\frac{\partial\mathbf{q}( \pmb{\zeta}_{2,j})}{\partial \zeta_{2,1}}\\
\frac{3h}{12}\frac{\partial\mathbf{q}( \pmb{\zeta}_{1,j})}{\partial \pmb{\zeta}_{1,j}} & 
\frac{5h}{12}\frac{\partial\mathbf{q}( \pmb{\zeta}_{2,j})}{\partial \pmb{\zeta}_{2,j}} \end{bmatrix}}_{A_{2}}\begin{bmatrix}\frac{\partial \pmb{\zeta}_{1,j}}{\partial \mathbf{x}_{j-1}}\\ \frac{\partial \pmb{\zeta}_{2,j}}{\partial \mathbf{x}_{j-1}} \end{bmatrix}.
\label{trapezoidalPD2}
\end{align}
Assuming that the matrix $(I_{2n}-A_{2})$ is invertible, where $A_{2}\in \mathbb{R}^{2n \times 2n}$, from the last expression we have 
\begin{align}
S =(I_{2n}-A_{2})^{-1} \widetilde{I}_n.
\label{trapezoidalPD3}
\end{align}
 After the matrix $S$ has been computed, we can substitute its elements in \eqref{trapezoidalPD1} to calculate the partial derivatives. 
\par
From \eqref{trapezoidalPD2} we see that to calculate the partial derivatives, we need to know the vectors $\{\pmb{\zeta}_{1,j}^{(i)},\pmb{\zeta}_{2,j}^{(i)} \}$, $j=0,1,\ldots, N-1$. These vectors can be obtained by simulating the system given by equations~\eqref{rk1},~\eqref{rk2}, and \eqref{rk3}, 
with an initial condition equal to $\mathbf{x}_{0}^{(i)}$.

\section{Optimal sensor selection}
\label{section4}

We consider a fixed number of sensors $r$ that 
represent directly observed nodes in the network. Accordingly, we introduce a vector $\mathbf{b}\in \{0, 1\}^n$, where its $i$th entry, denoted by $b_{i}$, is 
$1$ if node $i$ is observed and $0$ otherwise. In total, vector $\mathbf{b}$ should have $r$ entries equal to $1$, that is $\sum_{i=1}^{n} b_{i}= r$. 
Then, let the parametrized output equation be defined by $\mathbf{y}_{k}^{1}=C_{1}(\mathbf{b})\mathbf{x}_{k}$, where $C_{1}(\mathbf{b})=\text{diag}\left(\mathbf{b}\right)\in \mathbb{R}^{n\times n}$. Given a particular choice of $\mathbf{b}$, matrix $C_{1}$ is compressed into the matrix $C\in \mathbb{R}^{r\times n}$ by eliminating zero rows of $C_{1}$. Following the steps used to obtain equation~\eqref{nonlinearCompact3}, we define the parametrized equation $\mathbf{y}^{1}=\mathbf{w}^{1}\left(\mathbf{b},\mathbf{x}_{0}\right)$, where $\mathbf{w}^{1}\left(\mathbf{b},\mathbf{x}_{0}\right)=\text{col}\left(C_{1}(\mathbf{b})\mathbf{x}_{0},C_{1}(\mathbf{b})\mathbf{x}_{1},\ldots,C_{1}(\mathbf{b})\mathbf{x}_{N-1} \right)\in \mathbb{R}^{Nn\times n}$ and  $\mathbf{y}^{1}=\text{col}\left(\mathbf{y}_{0}^{1},\mathbf{y}_{1}^{1},\ldots, \mathbf{y}_{N-1}^{1} \right)\in \mathbb{R}^{Nn\times n}$. Linearizing the last equation around the initial state $\mathbf{x}_{0}$, we obtain
\begin{align}
\Delta \mathbf{y}^{1}=J^{1}(\mathbf{b},\mathbf{x}_{0})\Delta \mathbf{x_{0}},
\label{linearization1}
\end{align}
where $\Delta \mathbf{y}^{1} = \widetilde{\mathbf{y}} - \mathbf{y}^{1}\left(\mathbf{b},\mathbf{x}_{0}\right)$, $\widetilde{\mathbf{y}}$ is a vector of the output linearization,
$\mathbf{y}^{1}\left(\mathbf{b},\mathbf{x}_{0}\right)=\mathbf{w}\left(\mathbf{b},\mathbf{x}_{0}\right)$, 
$\Delta \mathbf{x}_{0}=\mathbf{x}_{0}^{*}-\mathbf{x}_{0}$, and
$\mathbf{x}^{*}_{0}$ is a state close to the initial state $\mathbf{x}_{0}$. Here, the 
Jacobian $J^{1}(\mathbf{b},\mathbf{x}_{0})\in \mathbb{R}^{Nn\times n}$ is calculated as 
\begin{align}
J^{1}\left(\mathbf{b},\mathbf{x}_{0}\right)&=\begin{bmatrix} C_{1}\left(\mathbf{b} \right) \\ C_{1}\left(\mathbf{b} \right)\frac{\partial \mathbf{x}_{1}}{\partial\mathbf{x}_{0}} \\ C_{1}\left(\mathbf{b} \right)\frac{\partial \mathbf{x}_{2}}{\partial \mathbf{x}_{1}}\frac{\partial \mathbf{x}_{1}}{\partial\mathbf{x}_{0}} \\ \vdots \\
C_{1}\left(\mathbf{b} \right)\frac{\partial \mathbf{x}_{N-1}}{\partial \mathbf{x}_{N-2}}\frac{\partial \mathbf{x}_{N-2}}{\partial \mathbf{x}_{N-3}}\ldots \frac{\partial \mathbf{x}_{1}}{\partial \mathbf{x}_{0}}  \end{bmatrix}=\left(I_{N}\otimes C_{1}\left(\mathbf{b}\right)\right)\underbrace{\begin{bmatrix} I_n\\ \frac{\partial \mathbf{x}_{1}}{\partial\mathbf{x}_{0}} \\ \frac{\partial \mathbf{x}_{2}}{\partial \mathbf{x}_{1}}\frac{\partial \mathbf{x}_{1}}{\partial\mathbf{x}_{0}} \\ \vdots \\
\frac{\partial \mathbf{x}_{N-1}}{\partial \mathbf{x}_{N-2}}\frac{\partial \mathbf{x}_{N-2}}{\partial \mathbf{x}_{N-3}}\ldots \frac{\partial \mathbf{x}_{1}}{\partial \mathbf{x}_{0}}  \end{bmatrix}}_{J^{2}\left(\mathbf{x}_{0}\right)},
\label{J1definition}
\end{align}
where $\otimes$ is the Kronecker product. This Jacobian has a form similar to the Jacobian given in equation~\eqref{jacobianOriginal}, except for the minus sign and the parametrization $\mathbf{b}$.
We note that the linearization \eqref{linearization1} can also be performed around states different from the initial state that we want to estimate. Importantly, our numerical experiments indicate that 
the approach is robust with respect to the uncertainties in the initial state.

The main idea of our approach is to quantify observability by the 
numerical ability to accurately solve
\eqref{nonlinearOptimization}, which depends to a large extent on the spectrum of $J^{1}\left(\mathbf{b},\mathbf{x}_{0}\right)^{T}J^{1}\left(\mathbf{b},\mathbf{x}_{0}\right)$. Our objective is to optimize this spectrum, which can be achieved, in particular, by maximizing the product of the eigenvalues.
Accordingly, we determine the optimal sensor selection as the solution of the following optimization problem:
\begin{align}
\begin{split}
& \min_{\mathbf{b}}\; \left\{  - \log \left[ \det \left(J^{1}\left(\mathbf{b},\mathbf{x}_{0}\right)^{T}J^{1}\left(\mathbf{b},\mathbf{x}_{0}\right) \right) \right]  \right\}  \\ 
& \text{subject to} \;\; \sum_{i=1}^{n} b_{i}= r, \; b_{i}\in \{0,1\},
\label{opt12}
\end{split}
\end{align}
where $r$ is the number of sensor nodes to be placed. Compared to other approaches in the literature based on 
empirical observability Gramians \cite{lall2002subspace,krener2009measures,powel2015empirical,qi2015optimal,singh2005determining} 
to quantify observability, our approach has the significant computational advantage that only a single simulation of the network dynamics is required to evaluate the cost function. In the case of empirical Gramians, typically order $n$ simulations of the dynamics are required
(for further details, see Appendix \ref{sec11}). 
We provide a thorough comparison of our approach for sensor placement with
other approaches
in Section \ref{section-results}.
It is also important to note that, in real systems, some variables cannot be experimentally observed at all. In our approach, such limitations can be easily incorporated as additional constraints in \eqref{opt12}.

Special attention in the objective function has to be given to the eigenvalues of the Jacobian. Specifically, one should always check a posteriori that the Jacobian does not have small or zero eigenvalues. If it does, a different observability measure should be used, such as the negative of the smallest eigenvalue or the condition number. A discussion of matrix measures other than the determinant is provided in Appendix \ref{sec11}.

The optimization problem in \eqref{opt12} belongs to the class of integer optimization problems and in the general case it is non-convex. To solve it, we use the NOMAD solver, implemented in the OPTI toolbox \cite{le2011algorithm,currie2012opti,qi2015optimal}. 
Our numerical results show that, despite being non-convex, the problem can be solved for $n=500$, $N=200$, and $r=50$ in just a few minutes on a standard current desktop computer.
We note that in larger systems the discrete optimization may become computationally infeasible. However, numerical solvers typically relax the problem by searching for an approximate (suboptimal) solution, which remains scalable. In practice, one can loosen the solver's convergence tolerances to obtain a fast approximate solution. In subsequent iterations, the solution may be further refined using tighter tolerances, until the solution (the set of chosen sensors) provides satisfactory state estimation results.

As numerically illustrated in the next section, the problem in \eqref{opt12} can also be effectively solved using the greedy algorithm proposed in \cite{summers2016submodularity}. However, the complete justification for using such an approach in the case of this optimization requires further theoretical developments that go beyond the scope of this paper (for more details, see Section \ref{sec12}).

\section{Numerical results}
\label{section-results}
We now demonstrate the efficacy of our approach through numerical experiments on chemical and biological networks. The general setup of 
the simulations is as follows.
The observation data set $\mathcal{O}_{N}$ is obtained by integrating the continuous-time dynamics in equation~\eqref{continiousTimeGeneral} using the MATLAB solver $\text{ode15s}$, specially designed for stiff dynamics.

\subsection{Applications to combustion networks}

We consider a hydrogen combustion reaction network ($H_2/O_2$) consisting of $9$ species and $27$ reactions, and the natural gas combustion network GRI-Mech 3.0, which consists of $53$ species and $325$ reactions. The state-space model is formulated as $\dot{\mathbf{x}}(t)=\Gamma \mathbf{q}_{c}\left(\mathbf{x}(t)\right)$, where $\mathbf{x}(t)\in \mathbb{R}^{n}_{\ge 0}$ is a state vector whose $n$ entries are the species concentrations. 
The matrix $\Gamma \in \mathbb{R}^{n\times n_{r}}$ consists of stoichiometric coefficients, where $n_{r}$ is the number of reactions, and the function $\mathbf{q}_{c}\left(\mathbf{x}(t)\right):\mathbb{R}^{n}\rightarrow \mathbb{R}^{n_r}$ consists of entries that are polynomials in $\mathbf{x}$ (see Appendix \ref{appendix-combust}).
The data to calculate the forward and backward reaction rates are taken from the reaction mechanisms database provided with the Cantera software \cite{Cantera}, files ``h2o2.cti'' and ``gri30.cti'' for the $H_{2}/O_{2}$ and GRI-Mech 3.0 networks, respectively. The rate constants are calculated using the modified Arrhenius law and the thermodynamic data available in the reaction mechanisms database. We ran the simulations under isothermal and constant volume conditions; the initial pressure used in our simulation is the atmospheric pressure and the temperature is $2500\,$K.  The entries of the ``true'' initial state $\mathbf{x}_{0}$  (concentrations of species in moles per liter) are calculated by
\begin{align}
\mathbf{x}_{0}=\mathbf{r}_{1}+\mathbf{c}_{1},
\label{generationx01}
\end{align}
where the entries of the vector $\mathbf{r}_{1}\in \mathbb{R}^{n}$ are drawn from the uniform distribution on the interval $(0,1)$ and $\mathbf{c}_{1}\in \mathbb{R}^{n}$ is the vector of ones.
The initial guesses $\mathbf{x}_{0}^{(0)}$ for solving the optimization problem in equation~\eqref{nonlinearOptimization} are generated in the same way (but independently from the true initial state). The lower bound  $\underline{\mathbf{x}}_{0}$ is a vector of zeros, whereas the upper bound $\overline{\mathbf{x}}_{0}$ is omitted in the optimization problem, because $\mathbf{x}_{0}$ is not bounded from above in this case.

\begin{figure}[h!]
\center
\includegraphics[width=4in]{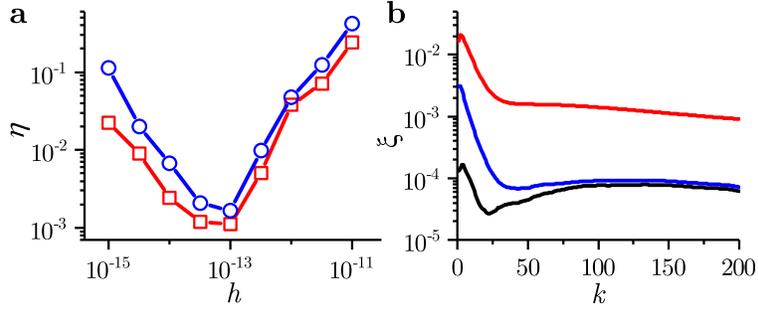}
\caption{\small Accuracy of the discretization methods. (\textbf{a}) Estimation error $\eta$ as a function of the discretization step size $h$, for the $H_{2}/O_{2}$ (squares) and GRI-Mech 3.0 (circles) networks, using the IRK model. (\textbf{b}) Comparison of the time-dependent error $\xi_{k} = \left\| \mathbf{x}_{k}-\mathbf{x}_{k}^{*} \right\|_{2} / \left\| \mathbf{x}_{k}^{*} \right\|_{2}$ for each step $k$, where $\mathbf{x}_{k} \in \mathbb{R}^{n}_{\ge 0}$ is computed using the BE (red line), TI (blue line), and IRK (black line) discretization methods. This comparison is for the $H_2/O_2$ network,         
with $N=200$, $h=10^{-13}$, and $\mathbf{x}_{k}^{*}$ 
computed using ode15s. In both panels, the results are averaged over $100$ realizations of random initial conditions, as defined in equation~\eqref{generationx01}.}
\label{figS2}
\end{figure}

The choice of the time step $h$ for discretizing the dynamics was determined based on its impact on the state estimation error, as shown in Fig.~\ref{figS2}. In general, the network cannot be \textit{fully} observed when the discretization time step is much larger than the dominant time constant in the system. We found that for the combustion networks analyzed here, the optimal choice is $h=10^{-13}$~(seconds). Since observing real combustion experiments at the resolution of picoseconds may be possible using femtosecond spectroscopy \cite{couris2014femtosec, li2013femtosec}, essentially all but the fastest reacting components could be estimated in a real experiment.

\begin{figure*}[h!]
\centering
\includegraphics[width=\textwidth]{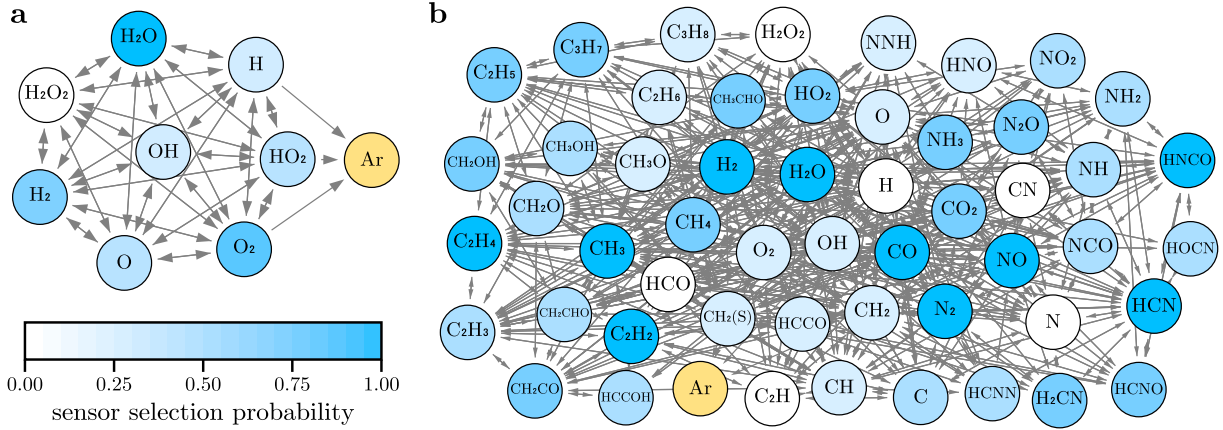}
\caption{\small OIDs and sensor selection probabilities for combustion networks. OIDs of the  (\textbf{a})  $H_{2}/O_{2}$ network and (\textbf{b})  GRI-Mech 3.0  network, where self-loops are omitted for clarity. Color-coded is the probability of selecting the node using the optimal sensor selection method (see text). Each network consists of two SCCs, one formed by Ar (always taken as a sensor) and the other by the remaining nodes. The data were computed using the IRK model for $N=200$ and $h=10^{-13}$. }
\label{fig2}
\end{figure*}

The OIDs are shown in Figs.~\ref{fig2}a and b for the $H_{2}/O_{2}$ and GRI-Mech 3.0 networks, respectively. In each case, the OID consists of two SCCs, of which one has a single node. This node is argon (Ar), which is an inert gas whose concentration remains constant and can influence the concentration of other species through pressure. The large SCC has no incoming edges and is thus a root component in both networks.  From the 
theory proposed in \cite{liu2013observability}, the complete network observability can be ensured by placing a single sensor in the root SCC. Accordingly, in our numerical simulations we make sure that at least one sensor is placed in the large SCC. In addition, before solving the optimization problem \eqref{opt12}, we place a sensor on the node forming the small SCC. This way, we avoid scenarios in which the observability measure would numerically overflow due to a badly conditioned Jacobian in \eqref{opt12} (especially when the number of sensor nodes is constrained to be small). More generally, some sensors may be placed a priori in non-root SCCs of any given network (irrespective of the number of SCCs) to avoid similar scenarios. In the context of our approach, these sensor additions can only further improve observability.

\begin{figure}[h!]
\centering
\includegraphics[width=4in]{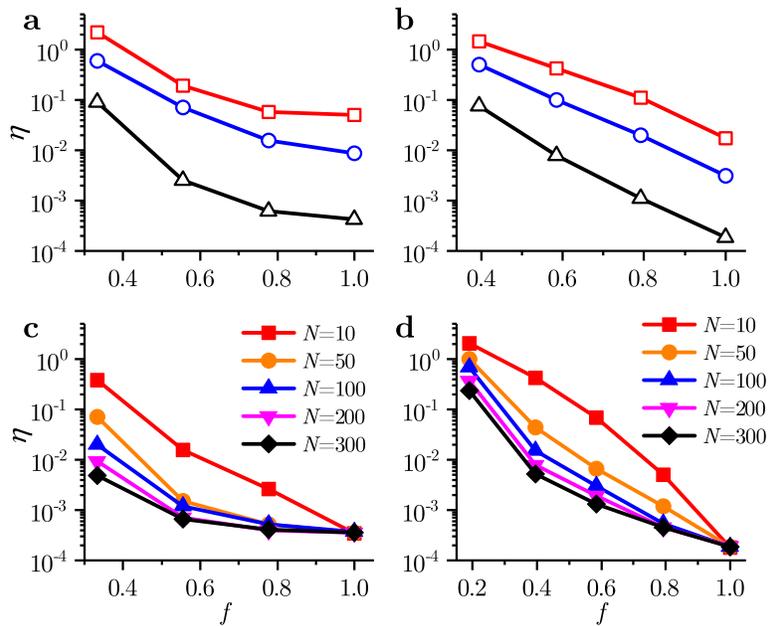}
\caption{\small Estimation error $\eta$ as function of the sensor fraction $f$ and observation length $N$. Results for the (\textbf{a}, \textbf{c}) $H_{2}/O_{2}$ network and (\textbf{b}, \textbf{d}) GRI-Mech 3.0 network.
 Panels
\textbf{a} and \textbf{b} compare the three models ($\vartriangle$-IRK, {\large $\circ$}-TI, and {\tiny $\square$}-BE) for $N=50$,
whereas  panels \textbf{c} and \textbf{d} compare different $N$ for the IRK model.
Each data point is an average over $100$ realizations of the random  sensor placement and 
initial guesses of the solution in the GRI-Mech 3.0 network and  over all possible placement configurations
in the $H_{2}/O_{2}$ network.  The discretization step was $h=10^{-13}$ in all simulations.
}
\label{fig3}
\end{figure}

To proceed, we define the sensor fraction $f = r/n$ as the fraction of
nodes in the network that  are selected as sensors.
We first compare the estimation errors $\eta$ for the three postulated models. The results are shown in Figs.~\ref{fig3}a and b for a random sensor selection in the large SCCs. 
These results show that the IRK model produces the lowest 
estimation error $\eta$, which is consistent with it being
more accurate than the BE and TI models. 
In addition, the IRK model also provides the smallest time-dependent error $\xi_k$ at every step $k$, as shown in Fig.~\ref{figS2}b.
The price of this is the increased computational complexity of the estimation procedure based on the IRK model. Note that all three models exhibit a non-zero estimation error for the sensor fraction $f=1$ (i.e., when all nodes are directly observed), which originates from model uncertainties, given that the postulated models only approximate the system used to generate the ``real'' data.

\begin{figure}[t!]
\centering
\includegraphics[width=4in]{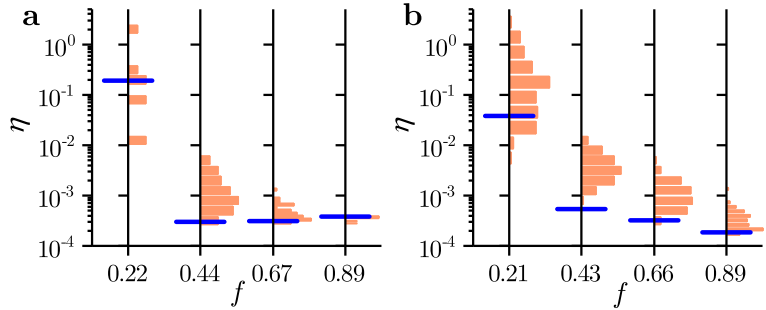}
\caption{\small Optimal sensor selection for the combustion networks. Estimation error $\eta$ for the (\textbf{a}) $H_{2}/O_{2}$ network and (\textbf{b}) GRI-Mech 3.0 network. For each network and sensor fraction $f$, the histogram 
presents $\eta$ for an exhaustive calculation in panel \textbf{a} and for
$100$ realizations of the random sensor placement 
in panel \textbf{b}. The blue line marks the estimation error for the calculated optimal sensor selection. The results are  generated using the IRK model for $N=200$ and $h=10^{-13}$.
}
\label{fig4}
\end{figure}

Next, we solve the optimal sensor placement problem given by equation~\eqref{opt12} 
for $N=200$ and several values of $f$, and  compute estimation error for the resulting optimal sensor selection. From the results shown in Figs.~\ref{fig4}a and b we conclude that the optimal sensor selection method performs significantly better than random sensor selection. 
Note that the observation horizon length $N$ has a strong influence on the sensor selection performance, which we show in detail in Section \ref{sec12}.

The probability of selecting the nodes using the optimal sensor selection method 
is obtained by calculating the frequency with which each node is chosen as a sensor in the solution of the optimal sensor selection problem for the various sensor fractions $f$ used in Fig.~\ref{fig4}. In this calculation, the contribution to the selection probabilities is weighted uniformly across different $f$, since the optimal sensor selection is generally unique for a given $f$. The probabilities are shown color-coded in Fig.~\ref{fig2}.

\subsection{Applications to biological networks}

\begin{figure*}[h!]
\centering
\includegraphics[width=\textwidth]{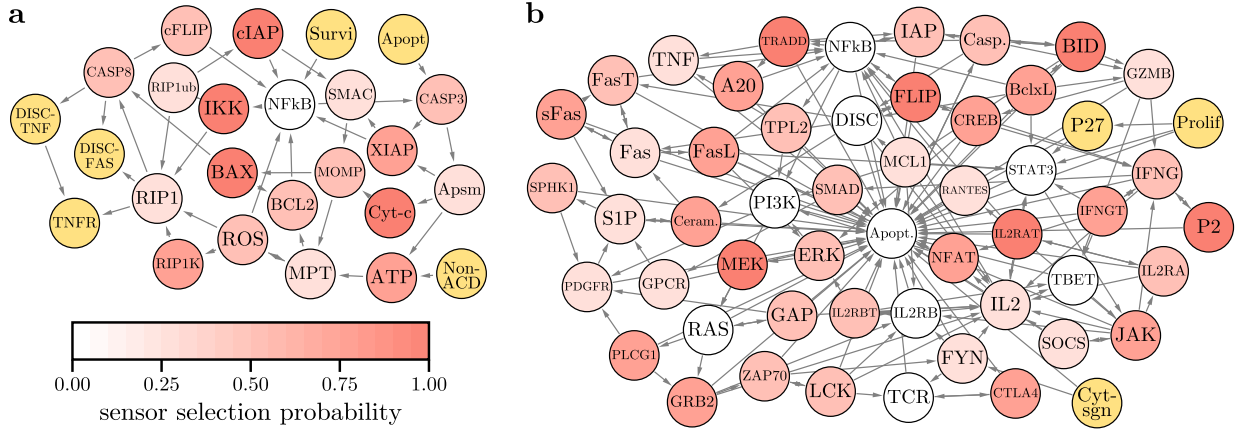}
\caption{\small OIDs and sensor selection probabilities for biological networks. 
Same as in Fig.~\ref{fig2} for the (\textbf{a}) CD network and (\textbf{b}) SS network. 
Yellow indicates single-node SCCs, whereas the remaining nodes belong to a giant SCC.
Sensors are placed on the yellow nodes a priori, which are then 
excluded from the optimal sensor selection.
The data were computed using the IRK model for $N=100$ and $h=0.02$ in panel \textbf{a}, and for $N=100$ and $h=0.05$ in panel \textbf{b}.}
\label{fig5}
\end{figure*}

We consider a cell death (CD) regulatory network model \cite{calzone2010mathematical} and a survival signaling (SS) network model of T cells
\cite{zhang2008network,saadatpour2011dynamical},
having $n=25$ and $n=54$ nodes, respectively.
Each node in these networks represents a gene, a gene product, or a concept (e.g., apoptosis), which can be active or inactive. Both systems are modeled as Boolean networks, where the activation of one node influences the activation of others. 
Following the standard practice, we transformed the Boolean relations into continuous-time dynamics 
using the $\text{ODEFY}$ software toolbox \cite{wittmann2009transforming}; specifically, we used the
Hillcube method with the threshold parameter of $0.5$. We also simplified the state-space model of the CD and SS networks, as follows. The original CD network contains $28$ nodes, however $3$ of them are input nodes (FASL, TNF, and FADD), whose values are set to $0.5$ in our simulations. This way, they are eliminated from the network. 
Similarly, the original SS network contains $60$ nodes, however $6$ of them are input nodes. The input nodes TAX, Stimuli2, and CD45 are set to zero, whereas the other input nodes, Stimuli, IL15, and PDGF, are set to $1$. 
After this simplification, the CD network has $7$ SCCs and the SS network has $4$ SCCs in the OID, as shown in Fig.~\ref{fig5}.

The entries of the initial state and the guesses for solving equation~\eqref{nonlinearOptimization} are generated from the uniform distribution on the interval $(0,1)$, where $\overline{\mathbf{x}}_{0}$ and $\underline{\mathbf{x}}_{0}$ are taken to be the vectors of ones and zeros, respectively.

\begin{figure}[h!]
\centering
\includegraphics[width=4in]{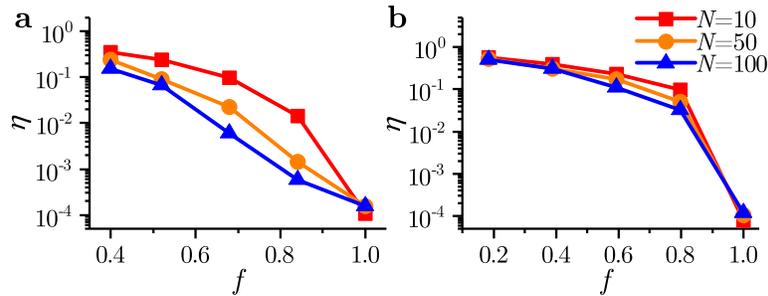}
\caption{\small Estimation error $\eta$ for the (\textbf{a}) CD network and (\textbf{b}) SS network as a function of the observation horizon and sensor fraction.  The results are generated for the IRK model with $h=0.02$ in panel \textbf{a} and $h=0.05$ in panel \textbf{b}. The results are averages over $100$ realizations of the random sensor selections and random initial conditions.}
\label{fig8signaling}
\end{figure}

\begin{figure}[h!]
\centering
\includegraphics[width=4in]{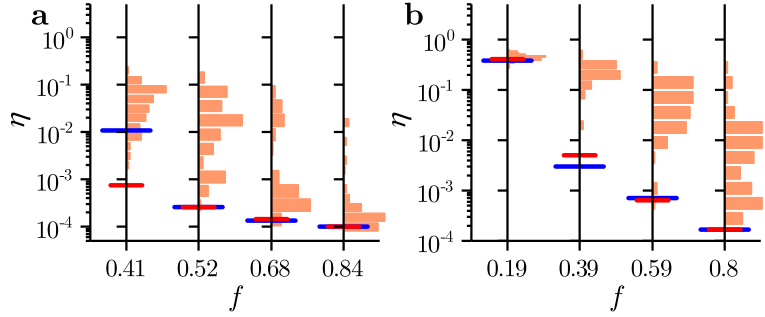}
\caption{\small Optimal sensor selection for the biological networks. Estimation error $\eta$ for the (\textbf{a}) CD network and (\textbf{b}) SS network. For each network and sensor fraction $f$, the histogram presents $\eta$ 
for $100$ realizations of the random sensor placement. 
The blue line marks the estimation error for the calculated optimal sensor selection; the red line indicates the corresponding result when information about the OID is ignored. 
The results are  generated using the IRK model for $N=100$ and $h=0.02$ in panel \textbf{a}, and for $N=100$ and $h=0.05$ in panel \textbf{b}.}
\label{fig4cd}
\end{figure}

The estimation error, shown in Fig.~\ref{fig8signaling}, behaves similarly to the error observed in the combustion networks. The results of the optimal sensor selection are shown in Figs.~\ref{fig4cd}a and b (blue lines), where we compare it with  two other approaches. As in Fig.~\ref{fig4}b, one approach is essentially graph-theoretic and consists of random placement under the constraint of having at least one sensor in each SCC
(histograms). The other approach is a variant of our optimal sensor selection method in which we lift the constraint of having at least one sensor in each SCC (red lines).  The results show the significant advantage of our approach compared to the others. The probability of selecting the nodes using the optimal sensor selection method is color-coded in Fig.~\ref{fig5}.

\subsection{Comparison to other sensor selection methods}
\label{sec12}

We now compare our sensor selection method to other methods
available in the literature. 
Two main approaches are available, one based on empirical Gramians (which can have various definitions) 
and the other based on heuristic measures of observability. 
Further details on these approaches are provided in
Appendix \ref{sec11}. For our analysis, we specify the following methods:
\begin{itemize}
\item Method 1. The empirical observability Gramian is computed using Definition 2 (see Appendix \ref{sec10}) and the optimal sensor selection is performed by solving the optimization problem in  equations~\eqref{opt1} and \eqref{opt2}.
\item Method 2. The empirical observability Gramian is computed using Definition 3 (see Appendix \ref{sec10}) and the optimal sensor selection is performed by solving the optimization problem in  equations~\eqref{opt1} and \eqref{opt2}.
\item Method 3. This is our approach. 
\item Method 4. Modified version of our approach, where we use a greedy algorithm \cite{summers2016submodularity} instead of the NOMAD method to solve \eqref{opt12}.
\end{itemize}
The motivation behind Method 4 is our conjecture that the cost function in \eqref{opt12} is submodular. If the conjecture is correct then the optimization can be solved 
by a greedy algorithm efficiently, and potentially faster than by the NOMAD solver. However, we leave the proof of this conjecture for future work.

Once the optimal sensor locations are determined, we estimate the initial state and compute the estimation error $\eta$. In order to validate the optimal sensor selection procedure, we compare it with a random sensor selection. Namely, in the case of the $H_{2}/O_{2}$ network, we generate all possible selections of the sensors for certain sensor fractions, whereas in the case of the GRI-Mech 3.0 network we generate $100$ random sensor selections, and for such selections we compute the initial state estimates and the estimation errors. In the validation step, the observation data 
are generated by simulating the network starting from the same initial state that has been used in sensor selection methods to compute the empirical  observability Gramians or the cost function in equation~\eqref{linObsMeasureG}. This is a ``true'' state, which is generated using equation~\eqref{generationx01}. The initial guess of the true state is also generated using equation~\eqref{generationx01}, and is generally not equal to the true state. We use $h=10^{-13}$, the IRK model, and vary the observation length $N$ to investigate the effect of the observation length on the optimal sensor selection performance. 
\par
\begin{figure}[h!]
\center
\includegraphics[width=4in]{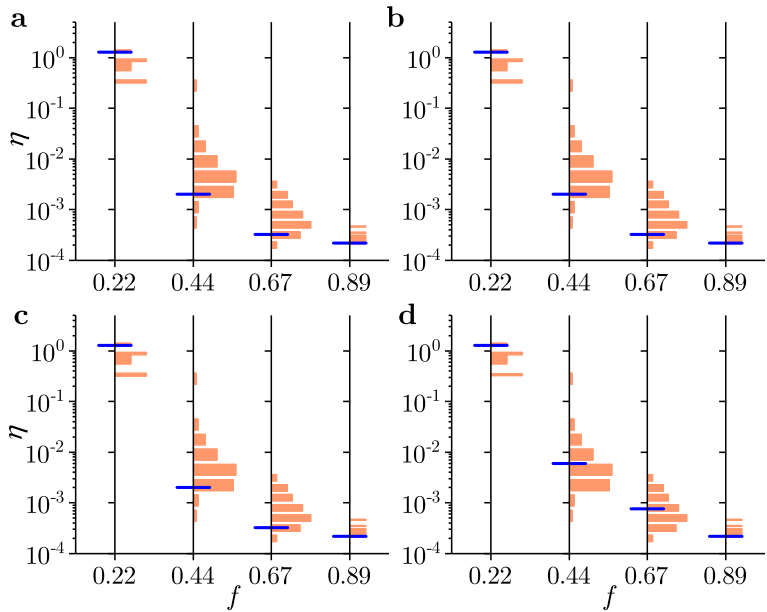}
\caption{\small  Four methods for sensor selection validated and compared on the $H_{2}/O_{2}$ network for short observation length. Estimation error $\eta$ for (\textbf{a}) Method 1, (\textbf{b}) Method 2, (\textbf{c}) Method 3, and (\textbf{d}) Method 4. The histograms correspond to the estimation errors for all possible combinations of the sensor nodes. The network is sufficiently small that exhaustive calculation of all combinations is possible in this case. The blue line represents the estimation error for the optimal sensor selection. The results are obtained for $N=50$, $h=10^{-13}$, and the IRK model.}
\label{figS6}
\end{figure}
\begin{figure}[h!]
\center
\includegraphics[width=4in]{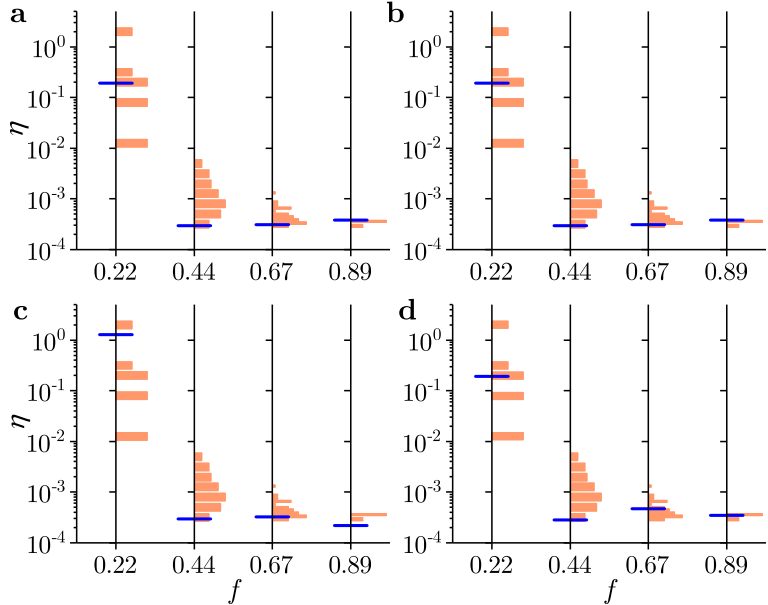}
\caption{\small Same as in Fig.~\ref{figS6} for the longer observation time of $N=200$.}
\label{figS7}
\end{figure}

Figure \ref{figS6} shows the comparison for the $H_{2}/O_{2}$ network and $N=50$. It can be observed that all $4$ methods perform relatively poorly compared to the random sensor selection. This is due to the fact that the total number of samples $N$ used in the computation of the cost functions and for the state estimation, is relatively short compared to the slowest time constant in the system. Consequently, the empirical observability Gramians and the Jacobians do not accurately capture the network dynamical behavior. Figure~\ref{figS7} shows the results for the $H_{2}/O_{2}$ network and $N=200$. In sharp contrast to the case of $N=50$, shown in Fig.~\ref{figS6}, the results shown in Fig.~\ref{figS7} are dramatically improved. We see that the methods perform well compared to a random sensor selection and that their performance is similar. The results can be additionally improved by selecting even larger $N$, however, at the expense of the increased computational complexity. To quantify the relative performance of the sensor selection methods, we calculate
the logarithmic error differences $\log \eta^{\mathrm{(method~}i\mathrm{)}} - \log \eta^{\mathrm{(method~3)}}$, for methods $i=1$, $2$, and $4$, respectively. These differences are $0.0$, $0.0$, and $-0.02$ for $f$$=$$0.44$. The corresponding differences are $-0.01$, $-0.01$, and $0.16$ for $f$$=$$0.67$. Finally, the differences are $0.24$, $0.24$, and $0.21$ for $f$$=$$0.89$.
It can be observed that, compared to the other Methods, Method 3 (our approach) 
has similar performance for $f$$=$$0.44$ and $f$$=$$0.67$, while it has a significantly better performance for $f$$=$$0.89$.

\subsubsection*{Comparison of computational complexity}
\begin{figure}[h!]
\center
\includegraphics[width=4in]{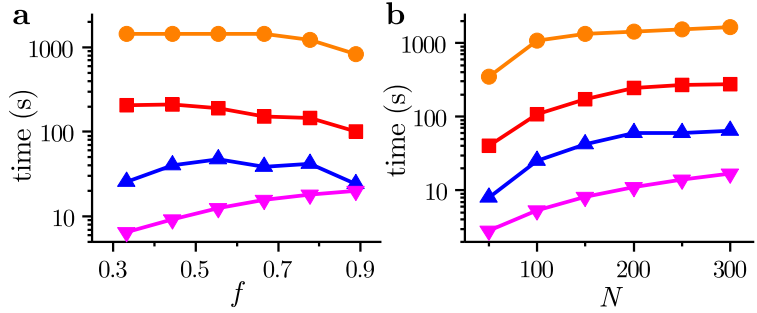}
\caption{\small Computational time for sensor selection as a function of (\textbf{a}) the fraction of observed nodes $f$ and (\textbf{b}) the observation horizon $N$.
The various curves correspond to Method 1 ({\tiny $\blacksquare$}), Method 2 ({\large $\bullet$}), Method 3 ($\blacktriangle$), and Method 4 ($\blacktriangledown$), applied to the GRI-Mech 3.0 network. Results are averaged over $10$ random initial guesses for the selected sensors, using the IRK model and $h=10^{-13}$, for $N=200$ in panel \textbf{a} and $f=0.5$ in panel \textbf{b}.}
\label{figS1}
\end{figure}

Compared to Methods 1 and 2, which are based on the observability empirical Gramians, the computational complexity of Methods 3 and 4 is much lower. This is because the computation of the empirical observability Gramians requires the network dynamics simulations for a number of perturbations of the initial condition, and each of these simulations have 
computational time scaling with at least $O(n^3)$. In practice, this computational complexity might be even higher due to model stiffness. On the other hand, to compute the cost functions for Methods 3 and 4, it is only necessary to simulate the dynamics for a single initial condition. 
This is reflected in the computational times of these methods shown in Fig.~\ref{figS1}. Method 3 (our approach) is almost an order of magnitude faster than Method 1, and two orders of magnitude faster than Method 2, across a wide range of parameters. We note that this difference will become even larger for larger networks. Method 4 (variant of our approach) is faster than Method 3 for the parameters shown in Fig.~\ref{figS1}, since the greedy algorithm 
evaluates the objective function fewer times than the NOMAD solver. However, this may not always be the case when both the network size $n$ and the fraction of observed nodes $f$ are large. In such cases, Method 4 may become slower because the number of objective function evaluations is proportional to $n \times f$ in this method. On the other hand, the number of objective function evaluations in Method 3 is determined by the convergence thresholds of the solver, which do not necessarily depend on $n$ or $f$.

\subsection{Discussion of the results}

\begin{figure}[h!]
\centering
\includegraphics[width=4in]{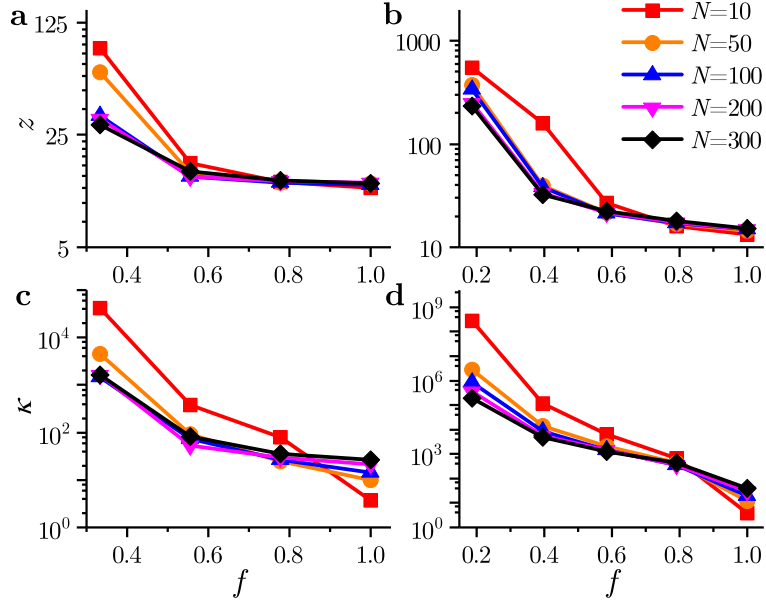}
\caption{\small Convergence and condition numbers of Jacobians of the nonlinear least-squares problem defined by equation~\eqref{nonlinearOptimization}. 
(\textbf{a,~b}) Number of iterations $z$ to compute the estimate as a function of the sensor fraction $f$ and the observation horizon $N$. (\textbf{c,~d}) Condition number $\kappa$ of the Jacobian matrix at the final estimate. Panels \textbf{a} and \textbf{c} correspond to the $H_{2}/O_{2}$ network, panels \textbf{b} and \textbf{d} correspond to the GRI-Mech 3.0 network. The parameters are the same as the ones used in Fig.~\ref{fig3}, panels \textbf{c} and \textbf{d}.}
\label{fig33}
\end{figure}

Figures~\ref{fig3}c and d indicate a trade-off between the sensor fraction $f$ and the observation length $N$ against the estimation error $\eta$.
The results show that there exists a fundamental obstacle in decreasing $\eta$ below a certain value for small $f$. Additional error reduction can be achieved by significantly increasing the observation length $N$, but our results show that $\eta$ starts to saturate for larger values of $N$. Furthermore, the value of $N$ is usually constrained by the  available computational power and the necessarily limited number of measurements in real experiments,
and it cannot be increased to infinity. On the other hand, for small $f$, the estimation procedure becomes ill-conditioned and the number of iterations of the TRR method dramatically increases (see Fig.~\ref{fig33}). This implies that, in practice, the observation performance can be severely degraded by measurement noise and that we need more time to compute the estimate when the fraction of observed nodes is smaller.

\begin{figure}[t!]
\center
\includegraphics[width=4in]{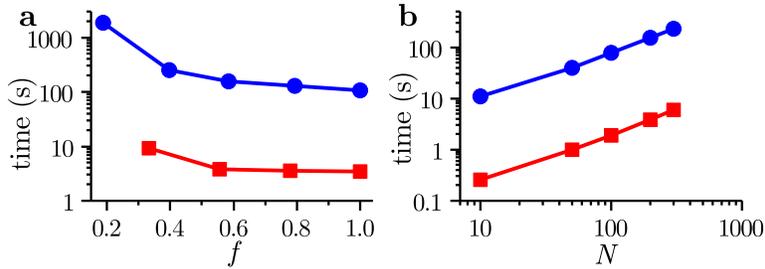}
\caption{\small Computational complexity of solving the nonlinear least-squares problem defined by equation~\eqref{nonlinearOptimization}. The computational complexity is shown for the $H_{2}/O_{2}$ ({\tiny $\blacksquare$}) and GRI-Mech 3.0 ({\large $\bullet$}) networks as a function of (\textbf{a}) the fraction of the observed nodes $f$ and of (\textbf{b}) the observation horizon $N$. In panel \textbf{a} the results correspond to an average of $100$ samples of randomly selected sensors and $N=200$, whereas in panel \textbf{b} the results correspond to $f=0.6$. In both panels, the results are obtained for the IRK model and $h=10^{-13}$.}
\label{figS5}
\end{figure}

For the networks we consider, our sensor selection procedure ensures that at least one sensor node is selected within each SCC of the OID. From a purely structural analysis of
Figs.~\ref{fig2} and \ref{fig5}, it would appear that 
the initial state could be reconstructed from the measurements of these sensors.
However, from our results, it follows that in practice this is not  
possible in general. As the sensor fraction $f$ approaches a small number, 
the estimation error $\eta$ dramatically increases. Furthermore, the number of iterations of the TRR method also dramatically increases, implying that the time to compute the estimate significantly increases (see Fig.~\ref{figS5}a).
We have verified through supplementary simulations that these conclusions are still valid even if the data-generating model is the same as the model postulated for observing the state. This suggests that the observed behavior does not originate from the model uncertainties alone---it can also be caused by limitations in numerical precision.
In the extreme case of a large network with a fraction of sensor nodes that is too small,
we would need an enormous amount of time to compute the estimate, and the estimation error would be large. In the case of small sensor fractions, the estimation error $\eta$ can be decreased by increasing $N$; however, larger $N$ also increases the computational complexity of the estimation method, as illustrated in Fig.~\ref{figS5}b.

For completeness, we have tested potential correlations with several node centrality metrics:
betweenness centrality, closeness centrality, degree (both in- and out-degrees), and pagerank,  which were calculated for every network studied here. In some cases we observe that the nodes selected by the optimal sensor selection method have larger than average centrality values,
 but beyond that we observe no clear correlation between individual centrality measures and the probability of selecting a node. The correlation coefficients between sensor selection probability and centrality measures are summarized in Appendix~\ref{section-centrality}. For the $H_{2}/O_{2}$ network, the most probable sensor nodes are the primary reactants and products, namely, $H_{2}$, $O_{2}$, and $H_{2}O$. For the GRI-Mech 3.0 network, this is only partially true; several probable sensor nodes are main reactants and products, such as $CH_{4}$, $O_{2}$, and $H_{2}O$, but others are unstable radicals, such as $CH_{3}$ and $C_{2}H_{5}$. We can only conclude that these species are selected because they carry the most information about the dynamics of other species in the mechanism. Similar conclusions hold for the biological networks as well. For example, in the CD network, two of the most frequently selected nodes (BAX and IKK) are primers for certain activation pathways \cite{calzone2010mathematical}.

\section{Conclusions}
\label{section-conclude}

Our results indicate that there is a fundamental limitation in network state estimation that, to the best of our knowledge, has not been considered before. This limitation stems from the complex interaction between three quantities: the number of available sensors, the observation length, and the condition number of the system's Jacobian. While in principle the estimation accuracy can be improved by increasing the number of sensors and the observation length, this is usually not feasible in practice. Apart from physical limitations of a real system (i.e., the number of available sensors, and the amount of available data), the fundamental limitation is the amount of computation needed to calculate the state estimate. Either the improvement of the condition number requires an unrealistic amount of data to be processed, or the number of iterations (and precision) remains too large due to bad conditioning. The best possible balance between these limitations for a given system can be found by our approach to optimal sensor selection. The results on the latter clearly illustrate the great potential of 
our framework compared to what would be possible with purely graph-theoretic 
approaches or approaches based on empirical Gramians, particularly as the number of SCCs and state variables increase. Our optimal sensor selection method indeed proves to be extremely useful in scenarios involving combinatorial explosion, such as in the case of  large networks.

\clearpage
\section*{Appendices}
\setcounter{section}{0}

\renewcommand\thesection{\Alph{section}}
\renewcommand\thesubsection{\thesection-\Alph{subsection}}

\section{Combustion networks}
\label{appendix-combust}

A combustion network is defined by 
\begin{align}
\sum_{i=1}^{n}\alpha_{ji}M_{i} \rightleftarrows    \sum_{i=1}^{n}\beta_{ji}M_{i}, \; j=1,2,\ldots, n_{r},
\label{chemicalReactionNetwork}
\end{align}
where $M_{i}$, $i=1,2,
\ldots, n$, are chemical species (e.g., $H$, $O$, or $H_{2}$), $n$ and $n_{r}$ are the total numbers of chemical species and reactions, respectively, and $\alpha_{ji}$ and $\beta_{ji}$ are stoichiometric coefficients. To equation~\eqref{chemicalReactionNetwork} we assign $n$ coupled differential equations of the form~\cite{perini2012analytical}
 \begin{align}
\dot{x}_{i}=\sum_{j=1}^{n_{r}}\left(\beta_{ji}-\alpha_{ji} \right)q_{j}\left(\mathbf{x}\right),\; i=1,2,\ldots, n, 
\label{stateSpaceChemicalReactions}
\end{align}
where $x_{i}\in \mathbb{R}_{\ge 0}$ is the concentration of the species $M_{i}$, $\mathbf{x}=\text{col}\left(x_{1},x_{2},\ldots, x_{n} \right)\in \mathbb{R}^{n}$ is the vector of concentrations (state vector), and $q_{j}$ is a polynomial function. This function has the following form:
\begin{align}
q_{j}\left(\mathbf{x}\right)=d_{j}^{(f)}\prod_{i=1}^{n}x_{i}^{\alpha_{ji}}-d_{j}^{(b)}\prod_{i=1}^{n}x_{i}^{\beta_{kj}},\;r=1,2,\ldots, n_{r},
\label{polynomialQ}
\end{align}
where $d_{j}^{(f)},d_{j}^{(b)} \in \mathbb{R}_{\ge 0}$ are the forward and backward reaction rates. The equations in \eqref{stateSpaceChemicalReactions} can be written compactly as
\begin{align}
\dot{\mathbf{x}}=\Gamma \mathbf{q}_{c}\left(\mathbf{x}\right),
\label{continiousTimeFinal}
\end{align}
where $\mathbf{q}_{c}\left(\mathbf{x}\right)=\text{col}\left(q_{1}\left(\mathbf{x} \right),q_{2}\left(\mathbf{x} \right),\ldots, q_{n_{r}}\left(\mathbf{x} \right) \right)$ and the matrix $\Gamma=[\gamma_{ij}]\in \mathbb{R}^{n\times n_{r}}$ has entries $\gamma_{ij}=\beta_{ji}-\alpha_{ji}$. Equation \eqref{continiousTimeFinal} is the state-space model of the combustion networks.


\section{Empirical Gramians of Nonlinear Systems}
\label{sec10}

Some of the approaches proposed in the literature for sensor 
placement are based on the concept of empirical observability Gramian of nonlinear systems  \cite{lall2002subspace,krener2009measures,powel2015empirical,qi2015optimal,singh2005determining}.
Here,
we review previous definitions of the empirical observability Gramian and introduce
a new definition that is more suitable for 
the class of network systems we consider. We also
provide guidelines for numerical computation and parameter selection.

\subsubsection*{Definition 1.}
We first consider the definition of the empirical observability Gramian introduced in \cite{lall2002subspace}. Let 
$\mathcal{T}^{n}=\{T_{1},T_{2},\ldots, T_{v} \}$ be a set of $v$ orthogonal, $n\times n$ matrices;
$\mathcal{M}=\{c_{1},c_{2},\ldots,c_{s} \}$ be a set of $s$ positive constants;
and let $\mathcal{E}^{n}=\{\mathbf{e}_{1},\mathbf{e}_{2},\ldots, \mathbf{e}_{n} \}$ be a set of $n$ standard unit vectors in $\mathbb{R}^{n}$.
Furthermore, let us introduce the mean $\overline{\mathbf{u}}$ of an arbitrary vector $\mathbf{u}$ as follows:
\begin{align}
\overline{\mathbf{u}}=\lim_{T\rightarrow\infty}\frac{1}{T}\int_{0}^{T}\mathbf{u}(t)\text{d}t.
\label{meanVectorU}
\end{align}
For the network dynamics with an output equation,
\begin{align}
\dot{\mathbf{x}}(t)&= \mathbf{q}\left(\mathbf{x}(t)\right),
\label{continiousTimeState} \\
\mathbf{y}(t)&=C\mathbf{x}(t), \label{continiousTimeOutput} 
\end{align}
the empirical observability Gramian  $\hat{X}_{1}\in \mathbb{R}^{n\times n}$ is defined by
\begin{align}
\hat{X}_{1}=\sum_{l=1}^{v}\sum_{m=1}^{s}\frac{1}{rsc_{m}^{2}}\int_{0}^{\infty}T_{l}\Psi^{lm}(t)T_{l}^{T}\text{d}t.\label{empiricalGramianDefinitionLall}
\end{align}
Here $\Psi^{lm}(t)\in \mathbb{R}^{n\times n}$ is a matrix whose $(i,j)$th entry is defined by
\begin{equation}
\Psi^{lm}_{i,j}(t)=\left(\mathbf{y}^{ilm}(t)-\overline{\mathbf{y}}^{ilm}\right)^{T}\left(\mathbf{y}^{jlm}(t)-\overline{\mathbf{y}}^{jlm}\right),
\label{entryPsi}
\end{equation} 
where $\mathbf{y}^{ilm}(t)$ is the output of the network corresponding to the initial condition $\mathbf{x}_{0}=c_{m}T_{l}\mathbf{e}_{i}$. The sets $\mathcal{M}$ and $\mathcal{T}^{n}$ are chosen by the user. Typical choices reported in the literature are \cite{qi2015optimal,lall2002subspace}: $\mathcal{M}=\{0.25,0.5,0.75,1\}$ and $\mathcal{T}^{n}=\{I_{n},-I_{n}\}$.

\subsubsection*{Definition 2.}
An alternative definition of the empirical observability Gramian \cite{powel2015empirical,krener2009measures} is
\begin{align}
\hat{X}_{2}(\tau)=\frac{1}{4\gamma^2}\int_{0}^{\tau}\Phi^{\gamma}(t)^{T}\Phi^{\gamma}(t) \text{d}t,
\label{gramianAlternative}
\end{align}
where 
\begin{align}
\Phi^{\gamma}(t)=\begin{bmatrix}\left( \mathbf{y}^{+1}(t)-\mathbf{y}^{-1}(t) \right) &  \ldots &\left( \mathbf{y}^{+n}(t)-\mathbf{y}^{-n}(t) \right) \end{bmatrix}. \label{phiDefinition}
\end{align}
Here the vector $\mathbf{y}^{\pm i}(t)$, $i=1,2,\ldots, n$, is the output of the network at the time $t$ for the initial condition $\mathbf{x}_{0}\pm \gamma \mathbf{e}_{i}$, where $\mathbf{x}_{0}$ is an arbitrary vector and $\gamma>0$ is a user choice. This definition is more attractive from the computational point of view, compared to the definition in equation~\eqref{empiricalGramianDefinitionLall}, because with the choice of the initial condition $\mathbf{x}_{0}$ we can freely choose the initial point around which we want to compute the empirical Gramian.

\subsubsection*{Definition 3.}
The initial conditions $\mathbf{x}_{0}\pm \gamma \mathbf{e}_{i}$ used to compute the empirical observability Gramian in equation~\eqref{gramianAlternative}, are not the most optimal choices from the computational point of view. Namely, $ \gamma \mathbf{e}_{i}$ perturbs $\mathbf{x}_{0}$ only in a single direction. We would like to construct a perturbation that is richer, such as the perturbation $c_{m}T_{l}\mathbf{e}_{i}$ used to compute the empirical observability Gramian in equation~\eqref{empiricalGramianDefinitionLall} (where $T_{l}$ is not an identity matrix). This motivates us to combine the above two definitions into a single one definition. We define the empirical observability Gramian as the following matrix:

\begin{align}
\hat{X}_{3}(\tau)=\sum_{l=1}^{v}\sum_{m=1}^{s}\frac{1}{4rsc_{m}^{2}} \int_{0}^{\tau}T_{l} \Phi^{lm}(t)^{T}\Phi^{lm}(t) T_{l}^{T} \text{d}t,
\label{definitionGramianNew}
\end{align}
where 
\begin{equation}
\Phi^{lm}(t)= \begin{bmatrix}\left( \mathbf{y}^{+1lm}(t)-\mathbf{y}^{-1lm}(t) \right),..., \left( \mathbf{y}^{+nlm}(t)-\mathbf{y}^{-nlm}(t) \right) \end{bmatrix},
\label{phiDefinition}
\end{equation}
and $\mathbf{y}^{\pm ilm}(t)$, $i=1,2,\ldots, n$, is the output of the network at time $t$ for the initial condition $\mathbf{x}_{0}\pm c_{m}T_{l} \mathbf{e}_{i}$, where $\mathbf{x}_{0}$ is an arbitrary vector.
\par
It can be easily shown that in the case of the linear dynamics
\begin{align}
\dot{\mathbf{x}}&=A\mathbf{x}, \label{stateLinear}  \\
\mathbf{y}&=C\mathbf{x}, \label{outputLinear}
\end{align}
where $A\in \mathbb{R}^{n\times n}$ and $C\in \mathbb{R}^{r\times n}$ are the constant system matrices, and for $\tau=\infty$ in \eqref{gramianAlternative}, the definitions in equations~\eqref{empiricalGramianDefinitionLall} and \eqref{gramianAlternative} become equal to the observability Gramian of the linear system given in equations~\eqref{stateLinear} and \eqref{outputLinear}:
\begin{align}
W=\int_{0}^{\infty} \exp(A^{T}t)C^{T}C\exp(At) \text{d}t,
\label{observabilityGramianLinearSystem}
\end{align}
where $\exp(At)\in \mathbb{R}^{n\times n}$ is the matrix exponential. In the sequel we show that the definition in equation~\eqref{definitionGramianNew} is equal to the definition given in equation~\eqref{observabilityGramianLinearSystem} for linear systems. For the linear dynamics in equations~\eqref{stateLinear} and \eqref{outputLinear}, and an arbitrary initial condition $\mathbf{z}\in \mathbb{R}^{n}$, we have
\begin{align}
\mathbf{y}(t)=C\exp(At)\mathbf{z}.
\label{initialStateResponse}
\end{align}
From equation~\eqref{initialStateResponse}, it follows that
\begin{align}
\mathbf{y}^{+ilm}(t)-\mathbf{y}^{-ilm}(t)&=C\exp(At)\left(\mathbf{x}_{0}+ c_{m}T_{l} \mathbf{e}_{i} \right)\notag \\& -C\exp(At)\left(\mathbf{x}_{0}-c_{m}T_{l} \mathbf{e}_{i} \right),\notag \\
&  =2c_{m}C\exp(At)T_{l}\mathbf{e}_{i}, 
\end{align}
which implies that the $(i,j)$th entry of the matrix $ \Phi^{lm}(t)^{T}\Phi^{lm}(t)$ in equation~\eqref{definitionGramianNew} is given by
\begin{align}
&\left( \mathbf{y}^{+ilm}(t)-\mathbf{y}^{-ilm}(t)\right)^{T}\left( \mathbf{y}^{+jlm}(t)-\mathbf{y}^{-jlm}(t)\right)\notag \\
&=4c_{m}^{2}\mathbf{e}_{i}^{T}T_{l}^{T}\exp(A^{T}t)C^{T}C\exp(At)T_{l}\mathbf{e}_{j}.
\label{derivation1tmp}
\end{align}
On the other hand, 
\begin{align}
\mathbf{e}_{i}^{T}Z\mathbf{e}_{j}=z_{i,j},
\label{derivation2tmp}
\end{align}
where $Z=[z_{i,j}]$ is an arbitrary matrix. From equations~\eqref{derivation1tmp} and \eqref{derivation2tmp}, we have
\begin{align}
T_{l} \Phi^{lm}(t)^{T}\Phi^{lm}(t) T_{l}^{T}=4c_{m}^{2}\exp(A^{T}t)C^{T}C\exp(At).
\label{derivation3tmp}
\end{align}
The last expression implies that
\begin{align}
\hat{X}_{3}(\tau)=&\sum_{l=1}^{v}\sum_{m=1}^{s}\frac{1}{rs} \int_{0}^{\tau}\exp(A^{T}t)C^{T}C\exp(At)   \text{d}t\notag \\
&=\int_{0}^{\tau}\exp(A^{T}t)C^{T}C\exp(At)   \text{d}t,
\label{definitionGramianNew22}
\end{align}
which for $\tau=\infty$ in equation~\eqref{definitionGramianNew} reduces to equation~\eqref{observabilityGramianLinearSystem}.

\subsubsection*{Computation of the observability Gramians}
To compute any of the empirical observability Gramians previously introduced, we first need to approximate the integrals. We use a trapezoidal method for approximating integrals because its computational complexity is low. To compute the expressions in equation~\eqref{definitionGramianNew}, we need to approximate the integral
\begin{align}
\chi=\int_{0}^{\tau}F^{lm}(t) \text{d}t,
\label{approximationIntegral}
\end{align}
where 
\begin{align}
F^{lm}(t)=T_{l} \Phi^{lm}(t)^{T}\Phi^{lm}(t) T_{l}^{T}.
\label{approximateIntegrals}
\end{align}
First we divide $[0,\tau]$ into $Q$ segments divided by points 
\begin{align}
0=t_{0}<t_{1}<t_{2}<\ldots <t_{Q}=\tau,
\label{pointsDivide}
\end{align}
The approximation of 
$\chi$
is defined by
\begin{align}
\chi \approx \frac{1}{2}\sum_{i=1}^{Q}\left(F^{lm}(t_{i-1})+F^{lm}(t_{i}) \right)\Delta t_{i},
\label{I1approximation}
\end{align}
where $\Delta t_{i}=t_{i}-t_{i-1}$. Taking this into account, the approximate empirical Gramian of $\overline{X}_{3}$, has the following form:
\begin{align}
\overline{X}_{3}=\sum_{l=1}^{v}\sum_{m=1}^{s}\frac{1}{8rsc_{m}^{2}} 
\sum_{i=1}^{Q}\left(F^{lm}(t_{i-1})+F^{lm}(t_{i}) \right)\Delta t_{i}.
\label{definitionGramianNewApprox}
\end{align}
Similarly, we can define the approximate Gramians of $\hat{X}_{1}$ and $\hat{X}_{2}$.
\par
To compute equation~\eqref{definitionGramianNewApprox}, we need to evaluate the matrix function 
$F^{lm}(t_{i})$
at the discrete-time steps $t_{i}$, $i=0,1,\ldots, Q$. This requires the knowledge of the state sequences $\mathbf{x}(t_{i})$ from different initial states. For simplicity, we choose equidistant time steps $\Delta t_{i}=\text{const}$. In the results reported in Section \ref{sec12}, the discrete-time sequences are computed on the basis of the IRK model.

\subsubsection*{Parameter choice of the observability Gramians}
To compute the empirical observability Gramians we need to choose the sets $\mathcal{T}^{n}$ and $\mathcal{M}$, the parameters $\gamma$ and $\tau$, as well as the initial condition $\mathbf{x}_{0}$. The principle for choosing these parameters is that all the initial states from which the computation of the state trajectory starts, should be within the physical limits of the state variables. For example, in the case of the combustion networks, the entries of the initial states should be positive, whereas in the case of the biological networks the entries of the initials states should be in the interval $[0,1]$.
For brevity, we explain the parameter selection for the case of combustion networks. Similar principles can be used in the case of biological networks. In all three definitions, the initial conditions are specified by equation~\eqref{generationx01}.

In Definition~1, we choose $\mathcal{T}^{n}=\{I_{n}\}$ (unlike the choice of $\{I_{n},-I_{n}\}$ reported in  \cite{qi2015optimal,lall2002subspace}) and $\mathcal{M}=\{0.25,0.5,0.75,1\}$. 
The improper integral is approximated by replacing the $\infty$ by a finite value $\tau$. A rule of thumb for the selection of the parameter $\tau$ is to choose it in such a way that for an arbitrary initial condition the majority of state trajectories approximately reach the steady states. However, because the computational complexity of computing the empirical observability Gramians increases with $\tau$, its value should not be very large. That is, there is a trade-off between the computational complexity and the value of $\tau$. As we show in Section \ref{sec12}, the value of $\tau$ dramatically influences the results, and it should be as large as possible. 
\par
In Definition 2, we choose $\gamma=0.5$. Our numerical results show that for the selection of the initial condition given by equation~\eqref{generationx01}, the entries of the perturbed initial conditions (almost) always stay positive. The value of $\tau$ is selected in the same manner as in Definition 1.
\par
In Definition 3, the matrices $T_{i}$ that are elements of the set $\mathcal{T}^{n}$ are chosen using the following procedure. First, we generate random matrices $S_{i}$, $i=1,2,\ldots, v$, whose entries are drawn from the standard normal distribution. After these matrices are constructed, we perform QR decompositions \cite{verhaegen2007}:
\begin{align}
S_{i}=Q_{i}R_{i}, \;\; i=1,2,\ldots, v.
\label{QRdecomposition}
\end{align}
The matrix $Q_{i}$ is orthogonal, and consequently, we chose $T_{i}=Q_{i}$ for all $i$.
With this choice of $T_{i}$ we are able to perturb more directions in the state-space around $\mathbf{x}_{0}$ compared to selecting $T_{i}$ as identity matrices. The set $\mathcal{M}$ is chosen as $\mathcal{M}=\{0.25, 0.5, 0.75, 1\}$. Our results show that for such a selection of $\mathcal{T}^{n}$ and $\mathcal{M}$, and for the selection of the initial condition given by equation~\eqref{generationx01}, the perturbed initial states almost always stay positive. The value of $\tau$ is selected in the same manner as in Definition 1.

\section{Existing approaches for optimal sensor selection}
\label{sec11}

We briefly summarize the two approaches from the literature against which our
optimal sensor selection method is compared in Section \ref{sec12}. 

The starting point of both approaches is the parameterized output equation: 
\begin{align}
\mathbf{y}_{k}^{1}=C_{1}(\mathbf{b})\mathbf{x}_{k},
\label{outputParameter}
\end{align}
where $C_{1}(\mathbf{b})=\text{diag}\left(\mathbf{b}\right)\in \mathbb{R}^{n\times n}$ and $\mathbf{b}\in \mathbb{R}^{n}$ is the parametrization vector. 
Once we have selected the sensors, the matrix $C_{1}$ is compressed into the matrix $C\in \mathbb{R}^{r\times n}$ that is used to define the output equation \eqref{observationEquation}. 

\subsubsection*{First approach}
The first approach is based on the empirical observability Gramians defined in Appendix \ref{sec10}. This approach is explained on the example 
of the empirical observability Gramian introduced in Definition 3. It can be easily modified such that it is based on the empirical Gramians introduced in the other two definitions.

We start from the approximate empirical Gramian $\overline{X}_{3}$ defined in equation~\eqref{definitionGramianNewApprox}. By substituting equation~\eqref{outputParameter} into equation~\eqref{phiDefinition}, the matrix  $\Phi^{lm}(t)$ becomes a function of the parametrization vector $\mathbf{b}$, that is, 
$\Phi^{lm}(t)=\Phi^{lm}\left(t,\mathbf{b}\right)$.
Then, substituting such an expression in equation~\eqref{approximateIntegrals}, we similarly obtain that $F^{lm}=F^{lm}(t,\mathbf{b})$.
This implies that the approximate Gramian in equation~\eqref{definitionGramianNewApprox} also becomes a function of the parametrization vector $\mathbf{b}$: 
$\overline{X}_{3}=\overline{X}_{3}\left(\mathbf{b}\right)$.

Several criteria have been used to quantify the degree of observability using empirical Gramians. Widely used criteria are the minimal singular (eigenvalue), trace, condition number, and the determinant of the empirical observability Gramian \cite{krener2009measures,qi2015optimal,singh2005determining}. The degree of observability based on the matrix determinant is attractive from the optimization point of view, mainly because the matrix determinant is a smooth function of the matrix entries
\cite{qi2015optimal}. Accordingly, we measure the degree of observability using the matrix determinant. Similarly to equation~\eqref{opt12}, 
the optimal sensor locations are determined as the solution of the following optimization problem:
\begin{align}
& \min_{\mathbf{b}}\; \{ - \log \left[ \det \left(\overline{X}_{3}\left(\mathbf{b} \right)  \right) \right] \} \label{opt1} \\
& \text{subject to} \;\; \sum_{i=1}^{n} b_{i}= r, \; b_{i}\in \{0,1\}. \label{opt2}
\end{align}
We solve this optimization problem using the NOMAD solver implemented in the OPTI toolbox \cite{le2011algorithm,currie2012opti,qi2015optimal}.

\subsubsection*{Second approach}
An elegant approach to the sensor selection for linear systems has been developed in \cite{summers2016submodularity}. This approach determines the optimal sensor locations by performing a finite number of evaluations of a set function that measures the observability degree of a network. It is shown that this method works well for linear systems provided that the set function is submodular (for more details, see \cite{summers2016submodularity}). Although it is developed for linear 
systems, we apply this approach to nonlinear systems without providing explicit proofs that justify the application of such a method for nonlinear systems. Our results show that this approach works relatively well even for nonlinear systems.
\par
Let $V=\{1,2,\ldots, n \}$ be a set where $n$ is the total number of nodes in our network. Let a set function $l:\; 2^{V}\rightarrow \mathbb{R}$ assign a real number to each subset $S$ of $V$. The subset $S$ is the set of sensors that we want to choose. For example, if we select $S=\{1,2,9\}$ then for such a choice, the function $l$ returns a real value quantifying the degree of observability. The problem of optimal sensor selection of $r$ nodes can be
 formulated as the following optimization problem:
\begin{align}
\max_{S\subset V,\; |S|=r} \; l\left(S \right).
\label{sensorSelectionOptimal}
\end{align}
It is easy to see that the optimization problems defined in equations~\eqref{opt1} and \eqref{opt2}, and in equations \eqref{opt12}, 
belong to the same class of optimization problems as the optimization problem in equation~\eqref{sensorSelectionOptimal}. The greedy algorithm for solving the optimization problem in equation~\eqref{sensorSelectionOptimal} has the following form. We start with an empty set $S_{0}$, then for $i=1,2,\ldots r$ and perform the following two steps:
\begin{enumerate}
\item Compute the gain $\Delta \left( a|S_{i} \right)=l\left(S_{i}\cup \{a\} \right)-l\left(S_{i} \right)$ for all elements $a\in V \setminus S_{i}$.
\item Define the set $S_{i+1}$ by adding the element $a$ to the set $S_{i}$ with the highest gain,
\begin{align}
S_{i+1} \leftarrow S_{i} \cup \{ \text{arg}\max_{a}\Delta \left( a|S_{i} \right)\;|\; a\in V\setminus S_{i} \}.
\label{highestGain}
\end{align}
\end{enumerate}
It is known that the greedy algorithm performs well for 
submodular functions $l$. In \cite{summers2016submodularity}, it is shown that the $\log\det\left(\cdot \right)$ function of the controllability Gramian of a linear system is submodular. Motivated by this we define the function $l$ to be
\begin{align}
l\left(\mathbf{b}\right)= - \log \left[ \det \left(J^{1}\left(\mathbf{b},\mathbf{x}_{0}\right)^{T}J^{1}\left(\mathbf{b},\mathbf{x}_{0}\right) \right) \right],
\label{linObsMeasureG}
\end{align}
where $\mathbf{x}_{0}$ is fixed and it is not considered as an argument of the $l$ function and $J^{1}\left(\mathbf{b},\mathbf{x}_{0}\right)$ is the Jacobian used in equation~\eqref{linearization1}. 
In equation~\eqref{linObsMeasureG} we slightly abuse the notation since $l$ is originally defined as a function that maps the set $S$ into a real value. On the other hand, the definition in equation~\eqref{linObsMeasureG} maps an $n$-dimensional vector into a real value. However, the sensors are marked by the position of the non-zero entries of $\mathbf{b}$. The proof of the submodularity of the function $l$ defined in equation~\eqref{linObsMeasureG} is left out and it is a future research topic.

\section{Sensor selection and centrality measures}
\label{section-centrality}

\begin{table}[H]
\caption{Correlation coefficients between the probability of selecting a node as a sensor (shown color-coded in Figs.~\ref{fig2} and \ref{fig5}) and various node's centrality measures in the OID.}
\label{table-centrality}
\centering
\begin{tabular*}{3.55in}{l r r r r}
\toprule
Centrality & $H_2/O_2$ & GRI-Mech 3.0 & CD & SS\\
\midrule
in-degree & 	$-0.57$ &	$-0.02$ &	$-0.63$	&	$-0.37$\\
out-degree & 	$-0.55$ &	$-0.04$ &	$-0.50$	&	$0.06$ \\
in-closeness & 	$-0.57$ &	$-0.05$ &	$-0.52$	&	$-0.40$ \\
out-closeness &	$-0.55$ &	$-0.19$ &	$-0.43$	&	$0.05$ \\
betweenness &	$-0.09$ &	$-0.13$ &	$-0.69$	&	$-0.38$ \\
pagerank &		$-0.57$ &	$-0.02$ &	$-0.03$	&	$-0.25$ \\
\bottomrule
\end{tabular*}
\end{table}



\vfill

\end{document}